# Activity Identification and Local Linear Convergence of Forward–Backward-type methods*

Jingwei Liang† Jalal M. Fadili† Gabriel Peyré‡

**Abstract.** In this paper, we consider a class of Forward–Backward (FB) splitting methods that includes several variants (*e.g.* inertial schemes, FISTA) for minimizing the sum of two proper convex and lower semi-continuous functions, one of which has a Lipschitz continuous gradient, and the other is partly smooth relative to a smooth active manifold $\mathcal{M}$. We propose a unified framework, under which we show that, this class of FB-type algorithms (i) correctly identifies the active manifolds in a finite number of iterations (finite activity identification), and (ii) then enters a local linear convergence regime, which we characterize precisely in terms of the structure of the underlying active manifolds. For simpler problems involving polyhedral functions, we show finite termination. We also establish and explain why FISTA (with convergent sequences) locally oscillates and can be slower than FB. These results may have numerous applications including in signal/image processing, sparse recovery and machine learning. Indeed, the obtained results explain the typical behaviour that has been observed numerically for many problems in these fields such as the Lasso, the group Lasso, the fused Lasso and the nuclear norm minimization to name only a few.

## 1 Introduction

### 1.1 Non-smooth optimization

In various fields of science and engineering, such as signal/image processing, inverse problems and machine learning, many problems can be cast as solving a *structured composite non-smooth optimization problem* of the sum of two functions, which usually reads

$$\min_{x\in\mathbb{R}^n} \ \Phi(x) \overset{\text{def}}{=} F(x) + R(x), \qquad (\mathcal{P}_{\text{opt}})$$

where

**(H.1)** $R \in \Gamma_0(\mathbb{R}^n)$, the set of proper convex and lower semi-continuous (lsc) functions on $\mathbb{R}^n$;

**(H.2)** $F \in C^{1,1}(\mathbb{R}^n)$, and the gradient $\nabla F$ is $\frac{1}{\beta}$-Lipschitz continuous;

**(H.3)** $\text{Argmin}(\Phi) \neq \emptyset$, *i.e.* the set of minimizers is non-empty.

From now on, we suppose that assumptions (**H.1**)-(**H.3**) hold. Problem ($\mathcal{P}_{\text{opt}}$) is closely related to finding solutions of the *monotone inclusion problem*

$$\text{Find } x \in \mathbb{R}^n \text{ such that } 0 \in A(x) + B(x), \qquad (\mathcal{P}_{\text{inc}})$$

†Normandie University, ENSICAEN, UNICAEN, GREYC, E-mail: {Jingwei.Liang, Jalal.Fadili}@ensicaen.fr.

‡CNRS, Ceremade, Université Paris-Dauphine, E-mail: Gabriel.Peyre@ceremade.dauphine.fr.

where we have

**(H.4)** $A : \mathbb{R}^n \rightrightarrows \mathbb{R}^n$ is a set-valued maximal monotone operator (see (A.1));

**(H.5)** $B : \mathbb{R}^n \to \mathbb{R}^n$ is maximal monotone and $\beta$-cocoercive (see (A.2));

**(H.6)** $\mathrm{zer}(A + B) \neq \emptyset$, *i.e.* the set of zeros of $A + B$ is non-empty.

For problem ($\mathcal{P}_{\mathrm{opt}}$), given a global minimizer $x^\star \in \mathrm{Argmin}(\Phi)$, then the corresponding first-order optimality condition reads

$$0 \in \partial R(x^\star) + \nabla F(x^\star),$$

where $\partial R$ denotes the sub-differential of $R$ at $x^\star$ (see definition (1.5)). Clearly, if we let $A = \partial R$ and $B = \nabla F$, then ($\mathcal{P}_{\mathrm{opt}}$) is simply a special case of ($\mathcal{P}_{\mathrm{inc}}$).

In this paper, our main focus is the non-smooth optimization problem ($\mathcal{P}_{\mathrm{opt}}$). Though some of our results are also valid for the monotone inclusion problem ($\mathcal{P}_{\mathrm{inc}}$), for instance the proposed Algorithm 1 and its global convergence analysis, see Theorem 2.1 and 2.3 in Section 2.

### 1.2 Forward–Backward-type splitting methods

The Forward–Backward (FB) splitting method [40] is a powerful tool for solving optimization problems ($\mathcal{P}_{\mathrm{opt}}$) with the additively separable and "smooth + non-smooth" structure. The standard (non-relaxed) version of FB updates a new iterate $x_{k+1}$ based on the following rule, ($x_0 \in \mathbb{R}^n$ is chosen arbitrarily)

$$x_{k+1} \stackrel{\text{def}}{=} \mathrm{prox}_{\gamma_k R}\big(x_k - \gamma_k \nabla F(x_k)\big), \ \gamma_k \in [\underline{\epsilon}, 2\beta - \overline{\epsilon}], \tag{1.1}$$

where $\underline{\epsilon}, \overline{\epsilon} > 0$, and $\mathrm{prox}_{\gamma R}$ denotes the *proximity operator* of $R$ which is defined as

$$\mathrm{prox}_{\gamma R}(\cdot) \stackrel{\text{def}}{=} \min_{x \in \mathbb{R}^n} \frac{1}{2}\|x - \cdot\|^2 + \gamma R(x).$$

The scheme (1.1) recovers the gradient descent method when $R = 0$, and the classic Proximal Point Algorithm (PPA) [53] when $F = 0$. Global convergence of the sequence $(x_k)_{k\in\mathbb{N}}$ generated by FB method is well established in the literature, based on the property that the composed operator $\mathrm{prox}_{\gamma R}(\mathrm{Id} - \gamma\nabla F)$ is so-called averaged non-expansive [12]. Moreover, sub-linear $O(1/k)$ convergence rate of the sequence of objective values of FB is also established in *e.g.* [47, 16, 14].

**Inertial schemes and FISTA** In the literature, different variants of the FB method were studied, and a popular trend is the inertial schemes which aim to speed up the convergence property of FB. In [51], a two-step algorithm called the "heavy-ball with friction" method is studied for solving ($\mathcal{P}_{\mathrm{opt}}$) with $R = 0$. It can be seen as an explicit discretization of a nonlinear second-order dynamical system (oscillator with viscous damping). This dynamical approach to iterative methods in optimization has motivated increasing attention in recent years. For instance, in real Hilbert spaces, it is used in [4] for solving ($\mathcal{P}_{\mathrm{opt}}$) with $F = 0$ and [5] for solving ($\mathcal{P}_{\mathrm{inc}}$) with $B = 0$ yielding an intertial PPA method. The authors in [44, 8, 41] propose different inertial versions of the FB method for solving ($\mathcal{P}_{\mathrm{opt}}$) and/or ($\mathcal{P}_{\mathrm{inc}}$) in real Hilbert spaces.

On the other hand, in the context of convex optimization, the accelerated FISTA method was proposed in [14], based upon the seminal work of [45], which achieves $O(1/k^2)$ convergence rate for the sequence of objective functions. However, while iterates generated by the FB are convergent, the convergence of FISTA iterates has remained a long-standing open problem. This question was recently settled in [19], followed by [9] in the continuous dynamical system case. More precisely, for $\gamma_k \in ]0, \beta]$ and a sequence of inertial parameter that converges at an appropriate rate (*i.e.* in the Algorithm 1 below, set $a_k = b_k = \frac{k-1}{k+q}$, $q > 2$), these authors have established (weak in infinite-dimensional Hilbert spaces) convergence of the iterates sequence while maintaining the $O(1/k^2)$ rate on the objective values. This rate is actually even $o(1/k^2)$ as proved in [7].

**Algorithm 1:** A General Inertial Forward–Backward splitting

**Initial**: $\bar{a} \le 1$, $\bar{b} \le 1$, $\underline{\epsilon}, \bar{\epsilon} > 0$ such that $\underline{\epsilon} \le 2\beta - \bar{\epsilon}$. $x_0 \in \mathbb{R}^n$, $x_{-1} = x_0$.
**repeat**

Let $a_k \in [0, \bar{a}]$, $b_k \in [0, \bar{b}]$, $\gamma_k \in [\underline{\epsilon}, 2\beta - \bar{\epsilon}]$:

$$y_{a,k} = x_k + a_k(x_k - x_{k-1}), \;\; y_{b,k} = x_k + b_k(x_k - x_{k-1}), \tag{1.2}$$

$$x_{k+1} = \mathrm{prox}_{\gamma_k R}\big(y_{a,k} - \gamma_k \nabla F(y_{b,k})\big). \tag{1.3}$$

$k = k + 1$;

**until** *convergence*;

In this paper, we propose a generalized inertial Forward–Backward splitting method (iFB) which by form covers all the above existing inertial schemes as special cases, see Algorithm 1. More precisely, based on the choice of the inertial parameters $a_k$ and $b_k$, the proposed method recovers the following special cases:

- $a_k = 0$, $b_k = 0$: this is the original FB method [40];
- $a_k \in [0, \bar{a}], b_k = 0$: this is the case studied in [44] for ($\mathcal{P}_{\mathrm{inc}}$). In the context of optimization with $R = 0$, one recovers the heavy ball method with friction in [51];
- $a_k \in [0, \bar{a}]$, $b_k = a_k$: this corresponds to the work of [41] for solving ($\mathcal{P}_{\mathrm{inc}}$). If moreover restrict $\gamma_k \in ]0, \beta]$ and let $a_k \to 1$, then Algorithm 1 specializes to FISTA-type methods [14, 19, 9, 7] developed for optimization.

When $a_k, b_k$ satisfy $a_k \in [0, \bar{a}]$, $b_k \in ]0, \bar{b}]$, $a_k \neq b_k$, Algorithm 1 is new in the literature to the best of our knowledge.

**Remark 1.1.**

(i) Though Algorithm 1 is stated for the optimization problem ($\mathcal{P}_{\mathrm{opt}}$), it readily extends to solve the monotone inclusion problem ($\mathcal{P}_{\mathrm{inc}}$), for which step (1.3) reads

$$x_{k+1} \stackrel{\mathrm{def}}{=} J_{\gamma_k A}\big(y_{a,k} - \gamma_k B(y_{b,k})\big), \tag{1.4}$$

where $J_{\gamma A} \stackrel{\mathrm{def}}{=} (\mathrm{Id} + \gamma A)^{-1}$ denotes the resolvent of $\gamma A$.

(ii) Though they share the same form of iteration when $a_k = b_k$, a notable difference between the inertial schemes and FISTA method is the range of choice for the stepsize $\gamma_k$, which is $[\underline{\epsilon}, 2\beta - \bar{\epsilon}]$ for the inertial methods, while only $]0, \beta]$ can be afforded by FISTA. This may have some impact on the practical convergence of the algorithm, see Section 5.5 for more details.

For the rest of the paper, we use the terminology *FB-type methods* for any scheme in the form of Algorithm 1 such that sequence $(x_k)_{k \in \mathbb{N}}$ converges. This will encompass the inertial schemes (denoted iFB) that we propose, the original FB method of course, and the sequence convergent FISTA method [19, 9] that corresponds to the specific choice of the intertial sequences $a_k = b_k = \frac{k-1}{k+q}$, $q > 2$. It should be noted, however, that our global convergence analysis to be presented in Section 2 does not cover the case of FISTA, which requires a specific proof strategy as developed in [19, 9].

### 1.3 Contributions

The study of (local) linear convergence of FB-type methods in the absence of strong convexity has become an active field in recent years, see the related work below for details. In general, most of the existing work focus mainly on some special cases (*e.g.* $R = \|\cdot\|_1$ in ($\mathcal{P}_{\text{opt}}$)), and the proofs of the results heavily rely on the specific structure of the function $R$, which makes them rather difficult to extend to other cases. Therefore, it is important to present a unified analysis framework, and possibly with stronger claims. This is one of the main motivations of this work. To be more precise, this paper consists of the following contributions.

**A general class of intertial algorithms** We present a unified iFB splitting class of algorithms for solving ($\mathcal{P}_{\text{opt}}$). It can be viewed as a versatile explicit-implicit discretization of a nonlinear second-order dynamical system with viscous damping, and thus covers existing methods as special cases. We establish global convergence of the iterates, and also stability to errors.

**Finite activity identification** Under the additional assumption that function $R$ is partly smooth at $x^\star \in \mathrm{Argmin}(\Phi)$ relative to a $C^2$-smooth manifold $\mathcal{M}_{x^\star}$ (see Definition 3.1) and a *non-degeneracy condition* at $x^\star$, we show that any FB-type method to solve ($\mathcal{P}_{\text{opt}}$) has the *finite time activity identification property*. Meaning that, after a finite number of iterations, say $K$, the iterates $x_k \to x^\star$ built by the FB-type method belong to $\mathcal{M}_{x^\star}$ for all $k \geq K$.

**Local linear convergence** Exploiting this identification property, we then show that the FB-type methods, locally along the manifold $\mathcal{M}_{x^\star}$, exhibit a linear convergence regime. We characterize this regime and the corresponding rates precisely depending on the structure of the active manifold $\mathcal{M}_{x^\star}$. For instance, we provide sharp estimates for the convergence rate. If moreover problem ($\mathcal{P}_{\text{opt}}$) has the structure described in Section 5.2, where $F$ is quadratic and $R$ is polyhedral, then *finite termination* can be obtained.

For the sequence convergent FISTA method, we draw two major conclusions:

- Locally, FISTA can be *slower* than the FB method (*e.g.* see Figure 3);
- We provide an explanation of the local oscillatory behaviour of FISTA (*e.g.* see Figure 4);

we describe precisely how these situations occur. This gives an enlightening explanation of the usefulness of the so-called restart method to locally accelerate the convergence of FISTA used by many authors, for instance in sparse recovery [27, 48, 26]: the algorithm is restarted after a certain number of iterations (set more or less empirically), where the inertial sequence $a_k = b_k$ is reset to $0$. In our work, we establish exactly the oscillation period of the FISTA iteration.

Building upon our local linear convergence analysis, we provide some pratical acceleration procedures. Indeed, once finite identification happens, the non-smooth convex problem ($\mathcal{P}_{\text{opt}}$) becomes (locally) equivalent to a $C^2$ smooth problem in the (possibly non-convex) active manifold $\mathcal{M}_{x^\star}$. In turn, this opens the door to acceleration, especially to apply higher order methods such as Newton or non-linear conjugate gradient.

Several numerical results are reported that confirm all our theoretical findings.

### 1.4 Related work

Finite support identification and local linear convergence of FB for solving a special instance of ($\mathcal{P}_{\text{opt}}$) where $F$ is quadratic and $R$ the $\ell_1$-norm (so-called LASSO problem), though in infinite-dimensional setting, is established in [16] under either a restrictive injectivity assumption, or a non-degeneracy assumption which is a specialization of ours (see (ND)). A similar result is proved in [28], for $F$ being a smooth convex and locally $C^2$ function and $R$ the $\ell_1$-norm, under restricted injectivity and non-degeneracy assumptions. The $\ell_1$-norm is polyhedral, hence partly smooth function, and is therefore covered by our results. [3] proved local linear convergence of FB to solve ($\mathcal{P}_{\text{opt}}$) for $F$ satisfying restricted smoothness and strong convexity assumptions, and $R$ being a so-called convex decomposable regularizer. Again, the latter is a subclass of partly smooth functions, and their result is thus covered by ours. For example, our framework covers the total variation (TV) semi-norm and $\ell_\infty$-norm regularizers which are not decomposable. Local linear convergence rate of FB for nuclear norm regularization is studied in [33] under local strong convexity assumption. Local linear convergence of FISTA for the Lasso problem (*i.e.* ($\mathcal{P}_{\text{opt}}$) for $F$ quadratic and $R$ the $\ell_1$ norm) has been recently addressed, for instance in [58], and also [34] under some additional constraints on the inertial parameters. The proposed work is also a deeper and sharper extension of our previous result on FB [39].

In [30, 31, 29], the authors have shown finite identification of active manifolds associated to partly smooth functions for a few algorithms, namely the (sub)gradient projection method, Newton-like methods, the proximal point algorithm and the algorithm in [59]. Their work extends that of *e.g.* [63] on identifiable surfaces (see references therein for related work of Dunn, and Burke and Moré). The algorithmic framework we consider encompasses all the aforementioned methods as special cases. Moreover, in all these works, the local convergence behaviour was not studied.

### 1.5 Notations

Throughout the paper, $\mathbb{N}$ is the set of non-negative integers and $k \in \mathbb{N}$ is the index. $\mathbb{R}^n$ is the Euclidean space of $n$ dimension, and $\mathrm{Id}$ denotes the identity operator on $\mathbb{R}^n$. For a nonempty convex set $\Omega \subset \mathbb{R}^n$, $\mathrm{ri}(\Omega)$ and $\mathrm{rbd}(\Omega)$ denote its relative interior and boundary respectively, $\mathrm{aff}(\Omega)$ is its affine hull, and $\mathrm{par}(\Omega) = \mathbb{R}(\Omega - \Omega)$ is the subspace parallel to it. Denote $\iota_\Omega$ the indicator function of $\Omega$, $\sigma_\Omega$ its support function and $\mathrm{P}_\Omega$ the orthogonal projector onto $\Omega$. For a matrix $M$, $\ker(M)$ is its null-space. We also denote $\mathrm{P}_\Omega$ the orthogonal projector onto $\Omega$. For a linear operator $L : \mathbb{R}^m \to \mathbb{R}^n$, we denote $L_T = L \circ \mathrm{P}_T$, and $L^+$ its Moore-Penrose pseudo-inverse.

The sub-differential of a function $R \in \Gamma_0(\mathbb{R}^n)$ is the set-valued operator,

$$\partial R : \mathbb{R}^n \rightrightarrows \mathbb{R}^n,\, x \mapsto \big\{ g \in \mathbb{R}^n | R(y) \geq R(x) + \langle g,\, y - x \rangle,\, \forall y \in \mathbb{R}^n \big\}. \tag{1.5}$$

We denote

$$T_x \stackrel{\text{def}}{=} \mathrm{par}\big(\partial R(x)\big)^\perp. \tag{1.6}$$

**Paper organization** The rest of the paper is organized as follows. Global convergence of the proposed iFB method is presented in Section 2. Then in Section 3, we introduce the concept of partial smoothness, and prove the finite activity identification property of the FB-type methods. We then turn to local linear convergence analysis in Section 4. Some hints about acceleration are provided in Section 4.5, and numerical results on various popular examples are reported in Section 5.

## 2 Global convergence of the inertial Forward–Backward

In this section, we establish global convergence of the iterates provided by Algorithm 1. We will state our results (Theorem 2.1 and 2.3) for the finite dimensional optimization problem ($\mathcal{P}_{\text{opt}}$). In fact, our global convergence results can handle the more general monotone inclusion problem ($\mathcal{P}_{\text{inc}}$) in an infinite dimensional real Hilbert space, where *weak* convergence of the iterates sequence can be obtained. The proofs given in Section A are written for this general setting.

## 2.1 Exact case

**Theorem 2.1 (Conditional convergence).** *Suppose that Algorithm 1 is run with* $\bar{a} < 1$*, and sequences* $(a_k)_{k\in\mathbb{N}}, (b_k)_{k\in\mathbb{N}}$ *such that*

$$\textstyle\sum_{k\in\mathbb{N}} \max\{a_k, b_k\}\|x_k - x_{k-1}\|^2 < +\infty. \tag{2.1}$$

*Then, there exists* $x^\star \in \mathrm{Argmin}(\Phi)$ *such that the sequence* $(x_k)_{k\in\mathbb{N}}$ *of Algorithm 1 converges to* $x^\star$*.*

The proof of Theorem 2.1 in given in Section A.

**Remark 2.2.** If $\forall k \in \mathbb{N}$, $a_k \geq b_k$, then (2.1) reduces to the even simpler form

$$\textstyle\sum_{k\in\mathbb{N}} a_k\|x_k - x_{k-1}\|^2 < +\infty. \tag{2.2}$$

Note that this condition is also the one provided in [44, 41] to ensure global convergence.

The terminology "conditional convergence" used in Theorem 2.1 refers to the fact that for the convergence to occur, the sequences $(a_k)_{k\in\mathbb{N}}$ and $(b_k)_{k\in\mathbb{N}}$ can be chosen depending (conditionally) on $(x_k)_{k\in\mathbb{N}}$ in such a way that (2.1) holds. This can be enforced easily by a simple online updating rule such as, given $a \in [0, 1], b \in [0, 1]$,

$$a_k = \min\left\{a, c_{a,k}\right\},\ \ b_k = \min\left\{b, c_{b,k}\right\}, \tag{2.3}$$

where $c_{a,k}, c_{b,k} > 0$, and $\max\{c_{a,k}, c_{b,k}\}\|x_k - x_{k-1}\|^2$ is summable. For instance, one can choose $c_{a,k} = \frac{c_a}{k^{1+\delta}\|x_k - x_{k-1}\|^2}, c_a > 0, \delta > 0$ and similarly for $c_{b,k}$.

One can also devise choices of $(a_k)_{k\in\mathbb{N}}$ and $(b_k)_{k\in\mathbb{N}}$ that are independent of $(x_k)_{k\in\mathbb{N}}$, and still guarantee global convergence. We dub this *unconditional convergence*. The following result generalizes those in [5, 44, 41].

**Theorem 2.3 (Unconditional convergence).** *Let* $\gamma_k$*,* $a_k$ *and* $b_k$ *as in Algorithm 1. Assume that there exists a constant* $\tau > 0$ *such that either of the following holds,*

$$\begin{cases} (1 + a_k) - \frac{\gamma_k}{2\beta}(1 + b_k)^2 > \tau : a_k < \frac{\gamma_k}{2\beta}b_k, \\ (1 - 3a_k) - \frac{\gamma_k}{2\beta}(1 - b_k)^2 > \tau : b_k \leq a_k \ \ or \ \ \frac{\gamma_k}{2\beta}b_k \leq a_k < b_k, \end{cases} \tag{2.4}$$

*Then* $\sum_{k\in\mathbb{N}} \|x_k - x_{k-1}\|^2 < +\infty$*, and there exists* $x^\star \in \mathrm{Argmin}(\Phi)$ *such that the sequence* $(x_k)_{k\in\mathbb{N}}$ *of Algorithm 1 converges to* $x^\star$*.*

See Section A for the proof. Figure 1 shows graphically the conditions in Theorem 2.3. We let $\tau = 0.01$ and two different choices of $\gamma$ are considered. It can be observed that with $\gamma$ becoming bigger, the range of $a, b$ in (2.4) becomes smaller.

## 2.2 Stability to errors

We now discuss the stability of the iFB method to errors. More precisely, we consider the case where $\partial R(x)$ and $\nabla F(x)$ are computed approximately. Toward this goal, we recall a notion which is inspired by the $\varepsilon$-approximate sub-differential in convex analysis.

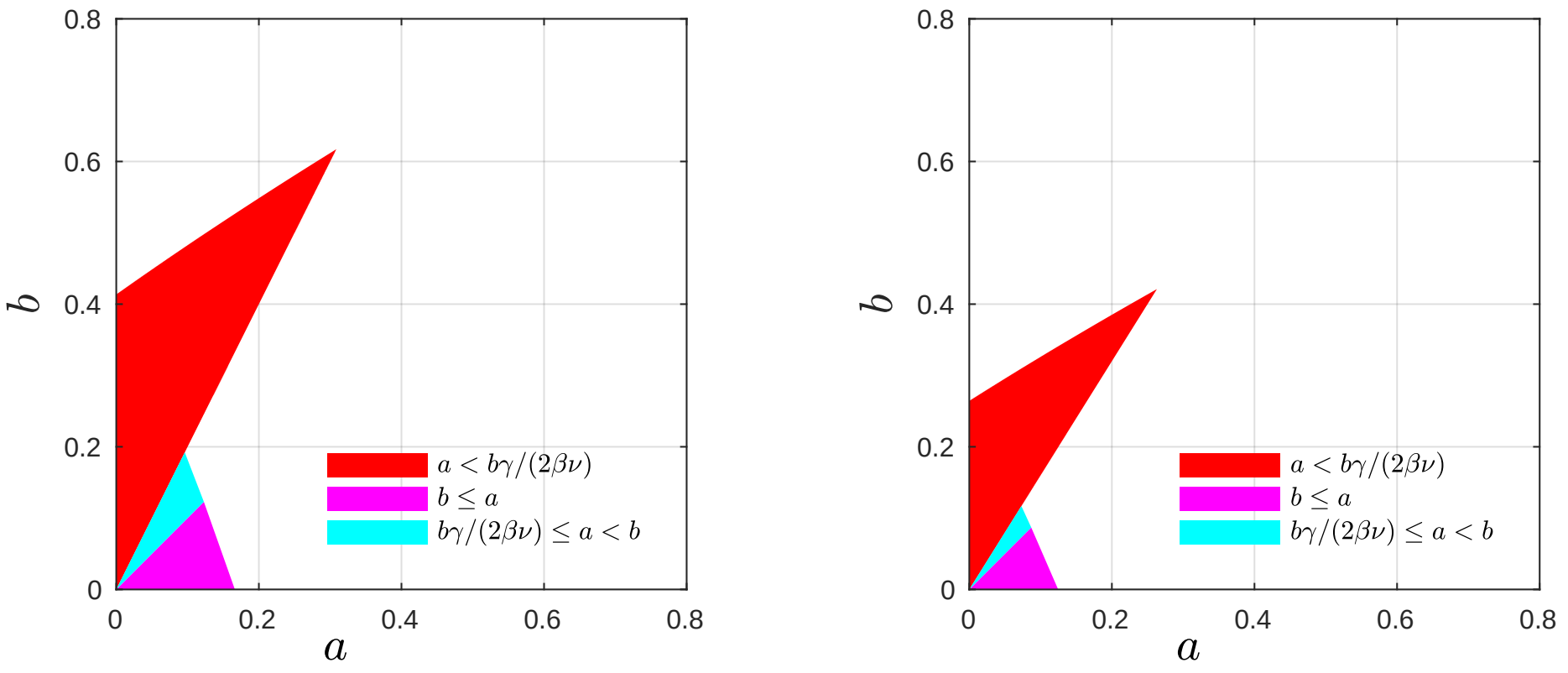

(a) $\gamma = \beta$ (b) $\gamma = 1.25\beta$

Figure 1: Sets of allowable $(a, b)$ ensuring the convergence for a given $\gamma$. (a) $\gamma = \beta$; (b) $\gamma = 1.25\beta$. We set the value of $\tau$ in (2.4) as $0.01$. Each color shaded region corresponds to a different condition appearing in (2.4), *i.e.* the cyan one corresponds to the first inequality of (2.4), while the magenta and red ones correspond to the two conditions of the second inequality of (2.4) respectively.

**Definition 2.4** (**$\varepsilon$-enlargement**)**.** Let $A : \mathbb{R}^n \rightrightarrows \mathbb{R}^n$ be a set-valued maximal monotone operator, $\varepsilon \geq 0$. Then the $\varepsilon$-enlargement of $A$ is defined as,

$$A^{\varepsilon}(x) \stackrel{\text{def}}{=} \big\{ v \in \mathbb{R}^n, \, \langle u - v, \, y - x \rangle \geq -\varepsilon, \; \forall y \in \mathbb{R}^n, \, u \in A(y) \big\}.$$

From the definition, for $0 \leq \varepsilon_1 \leq \varepsilon_2$ we have $A^{\varepsilon_1}(x) \subset A^{\varepsilon_2}(x)$ and $A^0(x) = A(x)$. Thus $A^{\varepsilon}$ is an enlargement of $A$.

Denote $\partial^{\varepsilon} R$ the $\varepsilon$-enlargement of $\partial R$. We now consider an inexact form of the iFB algorithm where step (1.3) is replaced by the corresponding inexact form that consists in finding $x_{k+1}$ such that

$$y_{a,k} - \gamma_k(\nabla F(y_{b,k}) + \xi_k) - x_{k+1} \in \gamma_k \partial^{\varepsilon_k} R(x_{k+1}), \tag{2.5}$$

where $\xi_k \in \mathbb{R}^n$ is the error in the evaluation of the gradient operator $\nabla F$. Observe that since the $\varepsilon$-approximate subdifferential of a proper closed convex function is contained in the $\varepsilon$-enlargement of its sub-differential [17], our setting also handles the case of approximate sub-differentials.

**Proposition 2.5.** *Consider Algorithm 1 with the inexact iteration* (2.5)*. Suppose that the conditions in Theorem 2.1 hold, and moreover, that one of the following holds,*

(i) $a_k \in ]0, \bar{a}]$, $\sum_{k \in \mathbb{N}} \varepsilon_k < +\infty$ *and* $\sum_{k \in \mathbb{N}} k \|\xi_k\| < +\infty$*;*
(ii) $a_k \equiv 0$, $\sum_{k \in \mathbb{N}} \varepsilon_k < +\infty$ *and* $\sum_{k \in \mathbb{N}} \|\xi_k\| < +\infty$*.*

*Then the conclusion of Theorem 2.1 holds true.*

See Section A for the proof. This result generalizes that of [44] who considered the case $b_k \equiv 0$ and $\xi_k \equiv 0$. [10] also studied the inexact sequence convergent FISTA method, *i.e.* $a_k = b_k = \frac{k-1}{k+q}$, $q > 2$, with the same errors as ours.

# 3 Partial smoothness and finite time activity identification

## 3.1 Partial smoothness

From now on, besides assumption (**H.1**), we assume that $R$ in ($\mathcal{P}_{\mathrm{opt}}$) is moreover *partly smooth function* relative to a smooth manifold. The notion of partial smoothness is first introduced in [37]. This concept, as well as that of identifiable surfaces [63], captures the essential features of the geometry of non-smoothness which are along the so-called active/identifiable manifold. For convex functions, a closely related idea is developed in [36]. Loosely speaking, a partly smooth function behaves smoothly as we move on the identifiable submanifold, and sharply if we move normal to the manifold. In fact, the behaviour of the function and of its minimizers depend essentially on its restriction to this manifold, hence offering a powerful framework for algorithmic and sensitivity analysis theory.

Let $\mathcal{M}$ be a $C^2$-smooth embedded submanifold of $\mathbb{R}^n$ around a point $x$. To lighten terminology, henceforth we shall state $C^2$-manifold instead of $C^2$-smooth embedded submanifold of $\mathbb{R}^n$. The natural embedding of a submanifold $\mathcal{M}$ into $\mathbb{R}^n$ permits to define a Riemannian structure on $\mathcal{M}$, and we simply say $\mathcal{M}$ is a Riemannian manifold. $\mathcal{T}_{\mathcal{M}}(x)$ denotes the tangent space to $\mathcal{M}$ at any point near $x$ in $\mathcal{M}$. More materials on manifolds are given in Section B.1.

We are now ready to state formally the class of partly smooth functions through its regularity properties.

**Definition 3.1 (Partly smooth function).** Let $R \in \Gamma_0(\mathbb{R}^n)$, $R$ is said to be *partly smooth at $x$ relative to a set $\mathcal{M}$* containing $x$ if $\partial R(x) \neq \emptyset$, and moreover

(i) **Smoothness**: $\mathcal{M}$ is a $C^2$-manifold around $x$, $R$ restricted to $\mathcal{M}$ is $C^2$ around $x$;

(ii) **Sharpness**: The tangent space $\mathcal{T}_{\mathcal{M}}(x)$ coincides with $T_x$ as given (1.6);

(iii) **Continuity**: The set-valued mapping $\partial R$ is continuous at $x$ relative to $\mathcal{M}$.

The class of partly smooth functions at $x$ relative to $\mathcal{M}$ is denoted as $\mathrm{PSF}_x(\mathcal{M})$.

One can easily show that a function in $\Gamma_0(\mathbb{R}^n)$ which is locally polyhedral around $x$ is partly smooth at $x$ relative to $x + T_x$. Polyhedrality also implies that the subdifferential is locally constant around $x$ along $x + T_x$. Capitalizing on the results of [37], it can be shown that under mild transversality conditions, the set of proper lsc convex and partly smooth functions is closed under addition and pre-composition by a linear operator. Moreover, absolutely permutation-invariant convex and partly smooth functions of the singular values of a real matrix, *i.e.* spectral functions, are convex and partly smooth spectral functions of the matrix [23]. Many examples of partly smooth functions that are popular in signal processing, machine learning and statistics will be discussed in Section 5.1.

[37, Proposition 2.10] allows to prove the following fact.

**Fact 3.2 (Local normal sharpness).** If $R \in \mathrm{PSF}_x(\mathcal{M})$, then all $x' \in \mathcal{M}$ near $x$ satisfy $\mathcal{T}_{\mathcal{M}}(x') = T_{x'}$. In particular, when $\mathcal{M}$ is affine or linear, then $T_{x'} = T_x$.

We now give expressions of the Riemannian gradient and Hessian (see Section B.1 for definitions) for the case of partly smooth functions relative to a $C^2$ submanifold. This is summarized in the following fact which follows by combining (B.2), (B.3), Definition 3.1, Fact 3.2 and [24, Proposition 17] (or [42, Lemma 2.4]).

**Fact 3.3.** If $R \in \mathrm{PSF}_x(\mathcal{M})$, then for any $x' \in \mathcal{M}$ near $x$

$$\nabla_{\mathcal{M}} R(x') = \mathrm{P}_{T_{x'}}(\partial R(x')),$$

and this does not depend on the smooth representation of $R$ on $\mathcal{M}$. In turn, for all $h \in T_{x'}$

$$\nabla^2_{\mathcal{M}} R(x')h = \mathrm{P}_{T_{x'}} \nabla^2 \widetilde{R}(x')h + \mathfrak{W}_{x'}\big(h, \mathrm{P}_{T_{x'}^{\perp}} \nabla \widetilde{R}(x')\big),$$

where $\widetilde{R}$ is a smooth extension (representative) of $R$ on $\mathcal{M}$, and $\mathfrak{W}_x(\cdot,\cdot) : T_x \times T_x^\perp \to T_x$ is the Weingarten map of $\mathcal{M}$ at $x$ (see Section B.1 for definitions).

### 3.2 Finite time activity identification

In this section, we state our result establishing that FB-type methods have the finite activity identification property.

**Theorem 3.4 (Finite activity identification).** *Suppose that the FB-type method is used to create a sequence $(x_k)_{k\in\mathbb{N}}$ that converges to $x^\star \in \mathrm{Argmin}(\Phi)$ such that $R \in \mathrm{PSF}_{x^\star}(\mathcal{M}_{x^\star})$, $F$ is locally $C^2$ around $x^\star$, and moreover the non-degeneracy condition*

$$-\nabla F(x^\star) \in \mathrm{ri}\big(\partial R(x^\star)\big), \tag{ND}$$

*holds. Then, there exists a large enough $K > 0$ such that for all $k \geq K$, $x_k \in \mathcal{M}_{x^\star}$.*

*If moreover,*

(i) *$\mathcal{M}_{x^\star}$ is an affine subspace, then $\mathcal{M}_{x^\star} = x^\star + T_{x^\star}$ and $y_{a,k}, y_{b,k} \in \mathcal{M}_{x^\star}, \forall k > K$;*

(ii) *$R$ is locally polyhedral around $x^\star$, then $y_{a,k}, y_{b,k} \in \mathcal{M}_{x^\star} = x^\star + T_{x^\star}$ for all $k > K$, $\nabla_{\mathcal{M}_{x^\star}} R(x_k) = \nabla_{\mathcal{M}_{x^\star}} R(x^\star)$, and $\nabla^2_{\mathcal{M}_{x^\star}} R(x_k) = 0$, $\forall k \geq K$.*

**Remark 3.5.**

(i) Recall that FB-type class of algorithms we consider contains the original FB method, the iFB one that we propose, and the FISTA method. The iFB is convergent under the assumptions of Theorem 2.1 or Theorem 2.3. The FISTA method is sequence convergent for $a_k = b_k = \frac{k-1}{k+q}$, $q > 2$, and $\gamma_k \equiv \gamma \in ]0, \beta]$; see [19, 9]. Thus, the finite identification property holds true for all these instances.

(ii) The non-degeneracy condition (ND) can be viewed as a geometric generalization of the strict complementarity of non-linear programming. Building on the arguments of [31], it is almost a necessary condition for the finite identification of $\mathcal{M}_{x^\star}$. Relaxing this assumption is a challenging problem in general.

**Proof.** Since $F$ locally is $C^2$ around $x^\star$, the smooth perturbation rule of partly smooth functions [37, Corollary 4.7], ensures that $\Phi \in \mathrm{PSF}_{x^\star}(\mathcal{M}_{x^\star})$.

By assumption, the sequence $(x_k)_{k\in\mathbb{N}}$ created by the FB-type method converges to $x^\star \in \mathrm{Argmin}(\Phi)$, and the latter is non-empty by assumption (**H.3**). Assumptions (**H.1**)-(**H.2**) entail that (ND) is equivalent to $0 \in \mathrm{ri}\big(\partial\big(\Phi(x^\star)\big)\big)$. Now (1.3) is equivalent to

$$\begin{aligned} & y_{a,k} - \gamma_k \nabla F(y_{b,k}) - x_{k+1} \in \gamma_k \partial R(x_{k+1}) \\ \iff & (a_k - b_k)(x_k - x_{k-1}) + (y_{b,k} - \gamma_k \nabla F(y_{b,k})) - (x_{k+1} - \gamma_k \nabla F(x_{k+1}) \in \gamma_k \partial\Phi(x_{k+1}). \end{aligned}$$

By Baillon-Haddad theorem [11], $\mathrm{Id} - \gamma_k\nabla F$ is averaged non-expansive for the prescribed range of $\gamma_k$, hence non-expansive, whence we get

$$\begin{aligned} \mathrm{dist}\big(0, \partial\Phi(x_{k+1})\big) &\leq \frac{1}{\gamma_k}\|(a_k - b_k)(x_k - x_{k-1}) + (y_{b,k} - \gamma_k\nabla F(y_{b,k})) - (x_{k+1} - \gamma_k \nabla F(x_{k+1}))\| \\ &\leq \frac{1}{\gamma_k}\big(|a_k - b_k|\|x_k - x_{k-1}\| + \|(y_{b,k} - \gamma_k\nabla F(y_{b,k})) - (x_{k+1} - \gamma_k\nabla F(x_{k+1}))\|\big) \\ &\leq \frac{1}{\gamma_k}\big(|a_k - b_k|\|x_k - x_{k-1}\| + \|x_k - x_{k+1}\| + b_k\|x_k - x_{k-1}\|\big) \\ &\leq \frac{1}{\gamma_k}\big(3\|x_k - x_{k-1}\| + \|x_k - x_{k+1}\|\big). \end{aligned}$$

Since $\liminf \gamma_k = \underline{\epsilon} > 0$ and $x_k$ is convergent, we obtain $\mathrm{dist}\big(0, \partial\Phi(x_{k+1})\big) \to 0$. Owing to assumptions (**H.1**)-(**H.2**), $\Phi$ is sub-differentially continuous at every point in its domain, and in particular at $x^\star$ for 0, which in turn entails $\Phi(x_k) \to \Phi(x^\star)$. Altogether, this shows that the conditions of [30, Theorem 5.3] are fulfilled, and the result follows.

(i) When the active manifold $\mathcal{M}_{x^\star}$ is an affine subspace, then $\mathcal{M}_{x^\star} = x^\star + T_{x^\star}$ owing to the normal sharpness property and the claim follows immediately;

(ii) When $R$ is locally polyhedral around $x^\star$, then $\mathcal{M}_{x^\star}$ is an affine subspace and the identification of $y_{a,k}, y_{b,k}$ follows from (i). For the rest, it is sufficient to observe that by polyhedrality, for any $x \in \mathcal{M}_{x^\star}$ near $x^\star$, $\partial R(x) = \partial R(x^\star)$. Therefore, combining Fact 3.2 and Fact 3.3, we get the second conclusion. ☐

**A bound on the identification iteration** In Theorem 3.4, we have not provided an estimate $K \geq 0$ beyond which finite identification occurs. There is of course a situation where the answer is trivial, *i.e.* $R$ is the indicator function of an affine subspace. However, knowing $K$ has practical interest, for instance, if one wants to switch to higher order acceleration (see Section 4.5). It is then legitimate to wonder whether such an estimate of $K$ can be given. In the following, we shall give a bound in some important cases. For the sake of simplicity, we state the result for the case of FB (*i.e.* $a_k = b_k \equiv 0$ in Algorithm 1). A similar reasoning can be easily generalized to the case of any converging FB-type method.

**Proposition 3.6.** *Suppose that the assumptions of Theorem 3.4 hold. Then the following holds.*

*(i) If the iterates are such that* $\partial R(x_k) \subset \mathrm{rbd}(\partial R(x^\star))$ *whenever* $x_k \notin \mathcal{M}_{x^\star}$*, then* $x_k \in \mathcal{M}_{x^\star}$ *for all* $k \geq \dfrac{\|x_0 - x^\star\|^2}{\underline{\epsilon}^2 \mathrm{dist}\big(-\nabla F(x^\star), \mathrm{rbd}(\partial R(x^\star))\big)^2}$;

*(ii) If $R$ is separable, i.e.* $R(x) = \sum_{i=1}^m \sigma_{C_i}(x_{b_i})$*, where* $\forall 1 \leq i \leq m, b_i \subset \{1,\ldots,n\}$, $\bigcup_{i=1}^m b_i = \{1,\ldots,n\}$*, and* $b_i \cap b_j = \emptyset, \forall i \neq j$*, and* $\dim(C_i) = |b_i|$*, then identification of* $\mathcal{M}_{x^\star}$ *occurs for some $k$ larger than* $\dfrac{\|x_0 - x^\star\|^2}{\underline{\epsilon}^2 \sum_{i \in I^c_{x^\star}} \mathrm{dist}\big(-\nabla F(x^\star)_{b_i}, \mathrm{rbd}(C_i)\big)^2}$*, where* $I_x \stackrel{\text{def}}{=} \big\{i : \; x_{b_i} \neq 0\big\}$.

**Proof.** (i) By firm non-expansiveness of $\mathrm{prox}_{\gamma_{k-1}R}$, and non-expansiveness of $\mathrm{Id} - \gamma_{k-1}\nabla F$, we have

$$
\begin{aligned}
\|x_k - x^\star\|^2 &\leq \|(\mathrm{Id} - \gamma_{k-1}\nabla F)(x_{k-1}) - (\mathrm{Id} - \gamma_{k-1}\nabla F)(x^\star)\|^2 \\
&\quad - \|x_{k-1} - \gamma_{k-1}\nabla F(x_{k-1}) - x_k + \gamma_{k-1}\nabla F(x^\star)\|^2 \\
&\leq \|x_{k-1} - x^\star\|^2 - \underline{\epsilon}^2\|u_k - \nabla F(x^\star)\|^2,
\end{aligned}
$$

where we denoted $u_k \stackrel{\text{def}}{=} (x_{k-1} - x_k)/\gamma_{k-1} - \nabla F(x_{k-1})$. By definition, we have $u_k \in \partial R(x_k)$. Suppose that identification has not occurred at $k$, *i.e.* that $x_k \notin \mathcal{M}_{x^\star}$, and hence $u_k \in \partial R(x_k) \subset \mathrm{rbd}(\partial R(x^\star))$. Therefore, continuing the above inequality, we get

$$
\begin{aligned}
\|x_k - x^\star\|^2 &\leq \|x_{k-1} - x^\star\|^2 - \underline{\epsilon}^2 \mathrm{dist}\big(-\nabla F(x^\star), \partial R(x_k)\big)^2 \\
&\leq \|x_{k-1} - x^\star\|^2 - \underline{\epsilon}^2 \mathrm{dist}\big(-\nabla F(x^\star), \mathrm{rbd}(\partial R(x^\star))\big)^2 \\
&\leq \|x_0 - x^\star\|^2 - k\underline{\epsilon}^2 \mathrm{dist}\big(-\nabla F(x^\star), \mathrm{rbd}(\partial R(x^\star))\big)^2,
\end{aligned}
$$

and $\mathrm{dist}(-\nabla F(x^\star), \mathrm{rbd}(\partial R(x^\star))) > 0$ owing to (ND). Taking $k$ as the largest integer such that the right hand is positive, we deduce that the number of iterations where identification has not occurred, does not exceed the given bound, whence our conclusion follows.

(ii) We have $\partial\sigma_{C_i}(x^\star_{b_i}) = C_i, \forall i \in I^c_{x^\star}$. In turn, by separability, $R$ is partly smooth at $x^\star$ relative to $\mathcal{M}_{x^\star} = \bigtimes_{i=1}^m \mathcal{M}_{x^\star_{b_i}}$, where $\mathcal{M}_{x^\star_{b_i}} = 0$ if $i \in I^c_{x^\star}$ and $\mathcal{M}_{x^\star_{b_i}} \neq 0$ otherwise. Suppose that at iteration $k$, $I^c_{x^\star} \cap I_{x_k} \neq \emptyset$. Denote $h_{k-1} = x_{k-1} - \gamma_{k-1}\nabla F(x_{k-1})$, and $h^\star = x^\star - \gamma_{k-1}\nabla F(x^\star)$. Thus for any $i \in I^c_{x^\star} \cap I_{x_k}$, we have

$$\begin{aligned} x_{k,b_i} - x^\star_{b_i} &= h_{k-1,b_i} - \mathrm{P}_{\gamma_{k-1}C_i}(h_{k-1,b_i}) \\ &= (h_{k-1,b_i} - h^\star_{b_i}) - (\mathrm{P}_{\gamma_{k-1}C_i}(h_{k-1,b_i}) - \mathrm{P}_{\gamma_{k-1}C_i}(h^\star_{b_i})) \end{aligned}$$

where we used Moreau identity in the first equality. Since $i \in I_{x_k} \cap I^c_{x^\star}$, we have $h_{k-1,b_i} \notin \gamma_{k-1}C_i$ and $h^\star_{b_i} \in \gamma_{k-1}C_i$, or equivalently, that $\mathrm{P}_{\gamma_{k-1}C_i}(h_{k-1,b_i}) \in \gamma_{k-1}\mathrm{rbd}(C_i) = \gamma_{k-1}\mathrm{rbd}(\partial\sigma_{C_i}(x^\star_{b_i}))$ and $\mathrm{P}_{\gamma_{k-1}C_i}(h^\star_{b_i}) = h^\star_{b_i}$. Combining this with the fact that the orthogonal projector on $\gamma_{k-1}C_i$ is firmly non-expansive, we obtain

$$\begin{aligned} \|x_{k,b_i} - x^\star_{b_i}\|^2 &\leq \|h_{k-1,b_i} - h^\star_{b_i}\|^2 - \|\mathrm{P}_{\gamma_{k-1}C_i}(h_{k-1,b_i}) - h^\star_{b_i}\|^2 \\ &= \|h_{k-1,b_i} - h^\star_{b_i}\|^2 - \|\mathrm{P}_{\gamma_{k-1}C_i}(h_{k-1,b_i}) + \gamma_{k-1}\nabla F(x^\star)_{b_i}\|^2 \\ &\leq \|h_{k-1,b_i} - h^\star_{b_i}\|^2 - \gamma^2_{k-1}\mathrm{dist}\big(-\nabla F(x^\star)_{b_i}, \mathrm{rbd}(C_i)\big)^2 \\ &\leq \|h_{k-1,b_i} - h^\star_{b_i}\|^2 - \underline{\epsilon}^2\mathrm{dist}\big(-\nabla F(x^\star)_{b_i}, \mathrm{rbd}(C_i)\big)^2. \end{aligned}$$

This bound together with non-expansiveness of $\mathrm{prox}_{\gamma_{k-1}C_i}$ and $\mathrm{Id} - \gamma_{k-1}\nabla F$ yield

$$\begin{aligned} \|x_k - x^\star\|^2 &= \textstyle\sum_{i\in I^c_{x^\star}} \|x_{k,b_i} - x^\star_{b_i}\|^2 + \sum_{j\in I_{x^\star}} \|x_{k,b_j} - x^\star_{b_j}\|^2 \\ &\leq \|h_{k-1} - h^\star\|^2 - \underline{\epsilon}^2 \textstyle\sum_{i\in I^c_{x^\star}} \mathrm{dist}\big(-\nabla F(x^\star)_{b_i}, \mathrm{rbd}(C_i)\big)^2 \\ &\leq \|x_{k-1} - x^\star\|^2 - \underline{\epsilon}^2 \textstyle\sum_{i\in I^c_{x^\star}} \mathrm{dist}\big(-\nabla F(x^\star)_{b_i}, \mathrm{rbd}(C_i)\big)^2 \\ &\leq \|x_0 - x^\star\|^2 - k\underline{\epsilon}^2 \textstyle\sum_{i\in I^c_{x^\star}} \mathrm{dist}\big(-\nabla F(x^\star)_{b_i}, \mathrm{rbd}(C_i)\big)^2, \end{aligned}$$

where the last term in the right hand side is strictly positive by (ND). Taking $k$ as the largest integer such that the right hand side is positive, we deduce that the number of iterations where $I^c_{x^\star} \cap I_{x_k} \neq \emptyset$ does not exceed the given bound. We then conclude that beyond this bound, there is no $i$ such that $\mathcal{M}_{x_{k,b_i}} \neq 0$ while $\mathcal{M}_{x^\star_{b_i}} = 0$. The proof is complete.

□

Note that, as intuitively expected, this bound increases as the non-degeneracy condition (ND) becomes more stringent. However, as it depends on $x^\star$, it is only of theoretical interest. In the separable case, observe that $\sum_{i\in I^c_{x^\star}} \mathrm{dist}\big(-\nabla F(x^\star)_{b_i}, \mathrm{rbd}(C_i)\big)^2 = \mathrm{dist}\big(-\nabla F(x^\star), \partial R(x^\star)\big)^2$ when $\sigma_{C_i}$ is differentiable at $x^\star_{b_i}$ for all $i \in I_{x^\star}$. The case of the $\ell_1$-norm considered in [28] is recovered in the second situation of Proposition 3.6 with $C_i \equiv [-\lambda, \lambda]$ for some $\lambda > 0$.

### 3.3 Stability to errors

Consider the inexact version (2.5) with $\varepsilon_k \equiv 0$, that is

$$x_{k+1} = \mathrm{prox}_{\gamma_k R}\big(y_{a,k} - \gamma_k(\nabla F(y_{b,k}) + \xi_k)\big).$$

Assume that $(\xi_k)_{k\in\mathbb{N}}$ is such that $(x_k)_{k\in\mathbb{N}}$ converges to some $x^\star \in \mathrm{Argmin}(\Phi)$ (see typically the summability conditions in Proposition 2.5(i)-(ii)). Then, since $\xi_k \to 0$, it can be easily seen from the proof of Theorem 3.4 that the activity identification property holds true for the above inexact iteration.

However, one cannot afford in general having non-zero errors $\varepsilon_k$ in the implicit step as in (2.5), even summable (see Proposition 2.5). The deep reason behind this is that in the exact case, under condition (ND), the proximal mappings of $R$ and $R + \iota_{\mathcal{M}_{x^\star}}$ locally agree nearby $x^\star$. This property is clearly violated if approximate proximal mappings are involved. Here is a simple example.

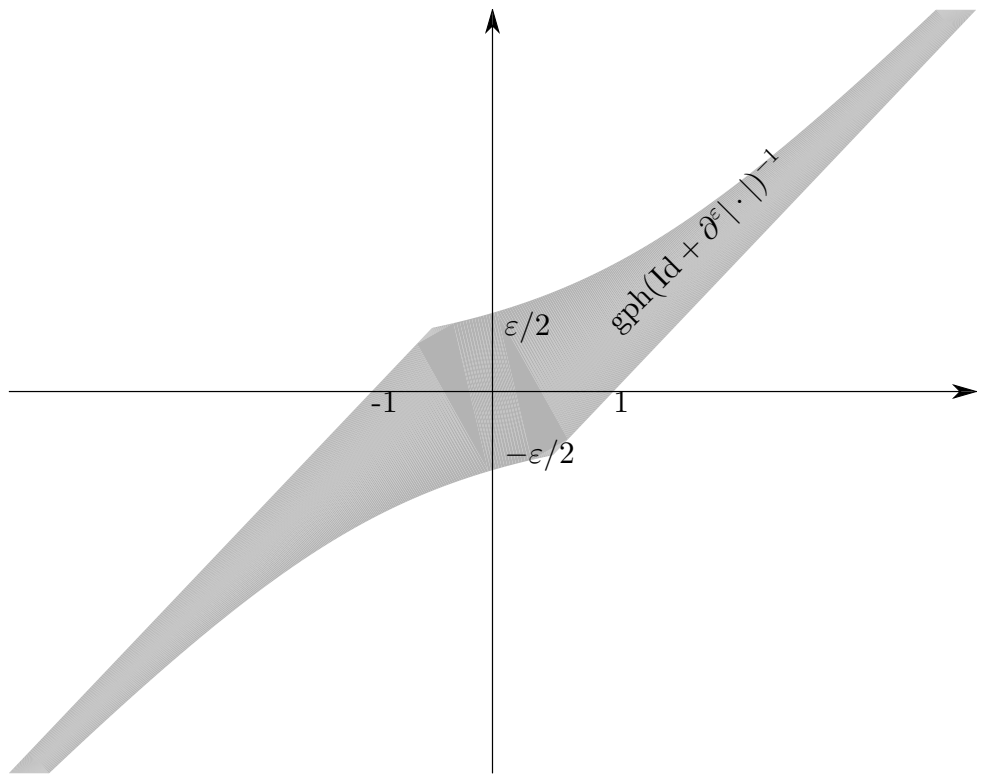

Figure 2: Graph of $(\mathrm{Id} + \partial^\varepsilon |\cdot|)^{-1}$.

**Example 3.7.** Let $F : x \in \mathbb{R} \mapsto \frac{1}{2}|\delta - x|^2$, with $\delta \in ]-1, 1[$, and $R : x \in \mathbb{R} \mapsto |x|$. It is easy to see that $\Phi \in \Gamma_0(\mathbb{R})$, and it has a unique minimizer $x^\star = \mathrm{prox}_{|\cdot|}(\delta) = \max(1 - 1/|\delta|, 0)\delta = 0$. Moreover, $\Phi$ is partly smooth at $x^\star$ relative $\mathcal{M}_{x^\star} = \{0\}$, and $\delta - x^\star = \delta \in \mathrm{ri}\big(\partial R(x^\star)\big) = ]-1, 1[$. Consider the inexact version of the FB algorithm

$$x_{k+1} \in (\mathrm{Id} + \partial^{\varepsilon_k} |\cdot|)^{-1}(\delta), \tag{3.1}$$

where we set $\gamma_k \equiv 1$, since $\nabla F$ is 1-Lipschitz. From [17, Example 5.2.5], we have

$$\partial^\varepsilon |\cdot|(x) = \begin{cases} [1 - \varepsilon/x, 1] & \text{if } x > \varepsilon/2 \\ [-1, 1] & \text{if } |x| \le \varepsilon/2 \\ [-1, -1 - \varepsilon/x] & \text{if } x < -\varepsilon/2, \end{cases}$$

whence the graph of $(\mathrm{Id} + \partial^\varepsilon |\cdot|)^{-1}$ can be easily deduced as displayed in Fig. 2. Thus, depending on $\varepsilon_k$ and the choice made in the inclusion (3.1), $x_k$ may never vanish for any finite $k$, *i.e.* $x_k \notin \mathcal{M}_{x^\star}$ for any finite $k$.

# 4 Local linear convergence of FB-type methods

We are now in position to present the local linear convergence result for FB-type methods, and all the proofs in this section are collected in Section B. Throughout this section, $x^\star$ is a global minimizer of problem ($\mathcal{P}_{\mathrm{opt}}$) such that the sequence $(x_k)_{k\in\mathbb{N}}$ provided by the FB-type method $x_k$ converges to $x^\star$. $\mathcal{M}_{x^\star}$ is the partial smoothness manifold of $R$ at $x^\star$, and $T_{x^\star}$ the corresponding tangent space.

**Restricted injectivity** In addition to the local $C^2$-smoothness assumption of $F$ made in Theorem 3.4, we suppose the following *restricted injectivity* condition,

$$\ker\big(\nabla^2 F(x^\star)\big) \cap T_{x^\star} = \{0\}. \tag{RI}$$

The local continuity of the Hessian of $F$ then implies that there exist $\alpha \geq 0$ and $\epsilon > 0$, such that $\forall h \in T_{x^\star}$,

$$\langle h,\, \nabla^2 F(x) h\rangle > \alpha \|h\|^2, \forall x \in \mathbb{B}_\epsilon(x^\star). \tag{4.1}$$

It turns out that under conditions (ND)-(RI), one can show that problem ($\mathcal{P}_{\text{opt}}$) admits a unique minimizer, and local quadratic growth of $\Phi$ if $R$ is moreover partly smooth. Recall that a function $\Phi$ grows quadratically locally around $x^\star$ if $\exists c > 0$ such that $\Phi(x) \geq \Phi(x^\star) + c\|x - x^\star\|^2$, $\forall x$ near $x^\star$.

**Proposition 4.1 (Uniqueness of the minimizer).** *Under assumptions* (**H.1**)-(**H.3**)*, let* $x^\star \in \mathrm{Argmin}(\Phi)$ *be a global minimizer of* ($\mathcal{P}_{\text{opt}}$) *such that* $F$ *is locally* $C^2$ *around* $x^\star$*. If conditions* (ND) *and* (RI) *are also fulfilled, then*

(i) $x^\star$ *is the unique minimizer of* ($\mathcal{P}_{\text{opt}}$).

(ii) *If moreover* $R \in \mathrm{PSF}_{x^\star}(\mathcal{M}_{x^\star})$*, then* $\Phi$ *has at least a quadratic growth near* $x^\star$.

**Remark 4.2.** In Proposition 4.1, partial smoothness of $R$ at $x^\star$ is not needed for the uniqueness claim (i). However, it brings more structure, hence the local quadratic growth property in (ii).

## 4.1 Locally linearized iteration

Define the following matrices which are all *symmetric*,

$$H \stackrel{\text{def}}{=} \gamma \mathrm{P}_{T_{x^\star}} \nabla^2 F(x^\star) \mathrm{P}_{T_{x^\star}}, \quad G \stackrel{\text{def}}{=} \mathrm{Id} - H, \quad U \stackrel{\text{def}}{=} \gamma \nabla^2_{\mathcal{M}_{x^\star}} \Phi(x^\star) \mathrm{P}_{T_{x^\star}} - H, \tag{4.2}$$

where $\nabla^2_{\mathcal{M}_{x^\star}} \Phi$ is the Riemannian Hessian of $\Phi$ on the manifold $\mathcal{M}_{x^\star}$ (see Fact 3.3).

**Lemma 4.3.** *For problem* ($\mathcal{P}_{\text{opt}}$)*, let* (**H.1**)-(**H.3**) *hold and* $x^\star \in \mathrm{Argmin}(\Phi)$ *such that* $R \in \mathrm{PSF}_{x^\star}(\mathcal{M}_{x^\star})$ *and* $F$ *is locally* $C^2$ *around* $x^\star$*. Then* $U$ *is symmetric positive semi-definite under either of the following circumstances:*

(i) (ND) *holds.*

(ii) $\mathcal{M}_{x^\star}$ *is an affine subspace.*

*In turn,* $\mathrm{Id} + U$ *is invertible, and* $W \stackrel{\text{def}}{=} (\mathrm{Id} + U)^{-1}$ *is symmetric positive definite with eigenvalues in* $]0, 1]$.

The following simple lemma gathers important properties of the matrices in (4.2).

**Lemma 4.4.** *For the matrices in* (4.2) *and* $W$*,*

(i) *Under* (**H.2**) *and* (RI)*,*

  (a) $H$ *is symmetric positive definite with eigenvalues in* $]\gamma\alpha, \frac{\gamma}{\beta}]$.

  (b) *For* $\gamma \in [\underline{\epsilon}, 2\beta - \overline{\epsilon}]$*,* $\underline{\epsilon}$ *and* $\overline{\epsilon} > 0$*,* $G$ *has eigenvalues in* $[-1 + \frac{\overline{\epsilon}}{\beta}, 1 - \alpha\underline{\epsilon}[ \subset ] - 1, 1[$.

  (c) *For* $\gamma \in [\underline{\epsilon}, \beta]$*,* $G$ *is also symmetric positive semi-definite with eigenvalues in* $[0, 1 - \alpha\underline{\epsilon}[ \subset [0, 1[$.

(ii) *If both the assumptions of Lemma 4.3 and (i) hold, then* $WG$ *has real eigenvalues lying in* $] - 1, 1[$*. If moreover* $\gamma \in [\underline{\epsilon}, \beta]$*, then* $WG$ *has eigenvalues lying in* $[0, 1[$.

Let $a \in [0, \bar{a}], b \in [0, \bar{b}], \gamma \in [\underline{\epsilon}, 2\beta - \overline{\epsilon}]$, define $r_k \stackrel{\text{def}}{=} x_k - x^\star, d_k \stackrel{\text{def}}{=} \begin{pmatrix} r_k \\ r_{k-1} \end{pmatrix}$, and the matrix

$$M \stackrel{\text{def}}{=} \begin{bmatrix} (a-b)W + (1+b)WG & -(a-b)W - bWG \\ \mathrm{Id} & 0 \end{bmatrix}. \tag{4.3}$$

Our interest in the vector $d_k$ is inspired by the convergence rate analysis of the heavy ball method [52, Section 3.2].

We now show that once the active manifold is identified, FB-type iteration locally linearizes.

**Proposition 4.5 (Locally linearized iteration).** *Let* (H.1)-(H.3) *hold, and assume that an FB-type method is used to create a sequence* $(x_k)_{k\in\mathbb{N}}$ *that converges to* $x^\star \in \mathrm{Argmin}(\Phi)$ *such that* (ND) *and* (RI) *hold. If moreover,*

$$a_k \to a \in [0,1],\ b_k \to b \in [0,1],\ \gamma_k \to \gamma \in [\underline{\epsilon}, 2\beta - \overline{\epsilon}], \tag{4.4}$$

*then for $k$ large enough, we have*

$$d_{k+1} = M d_k + o(\|d_k\|). \tag{4.5}$$

*The $o(\cdot)$ term disappears when $R$ is locally polyhedral and $(\gamma_k, a_k, b_k)$ are chosen constant.*

**Remark 4.6.**

(i) (4.4) asserts that both the inertial parameters $(a_k, b_k)$ and the step-size $\gamma_k$ should converge to some limit points, and this condition cannot be relaxed in general.

(ii) For the FB method (*i.e.* $a_k = b_k \equiv 0$), (4.3) can be further simplified, and the corresponding linearized iteration can be stated in terms of $r_k$ directly, which reads

$$r_{k+1} = WGr_k + o(\|r_k\|). \tag{4.6}$$

(iii) Proposition 4.5 also covers the sequence convergent FISTA method [19, 9], *i.e.* $a_k = b_k = \frac{k-1}{k+q}$, where $q > 2$ is a constant, and $\gamma_k \equiv \gamma \in ]0, \beta]$. In this case, we have indeed $a_k \to a = b = 1$.

## 4.2 Spectral properties of $M$

Our aim now is to establish local linear convergence of FB-type schemes. For this, given the structure of the locally linearized iteration (4.5), it is sufficient to strictly upper-bound by 1 the spectral radius of $M$, and conclude using standard arguments. This is what we are about to do.

The rationale is to start by relating explicitly the eigenvalues of $M$ to those of $G$ or $WG$, and then use Lemma 4.4 to upper-bound the spectral radius of $M$. However, given the structure of $M$, this is a challenging linear algebra problem, and can only be done for some cases: $a$ and $b$ possibly different but the the function $R$ is locally polyhedral, or $R$ is a general partly smooth function but $a = b$. These situations are not restrictive at all and cover all interesting applications we have in mind.

Let $\eta$ and $\sigma$ be an eigenvalue of $WG$ and $M$ respectively. We denote $\underline{\eta}, \overline{\eta}$ the smallest and largest (signed) eigenvalues of $WG$, and $\rho(M)$ the spectral radius of $M$.

**Locally polyhedral case** When $R$ is locally polyhedral, $U$ vanishes and $W = \mathrm{Id}$, then $M$ in (4.3) simplifies to the following form

$$M = \begin{bmatrix} (a-b)\mathrm{Id} + (1+b)G, & -(a-b)\mathrm{Id} - bG \\ \mathrm{Id}, & 0 \end{bmatrix}. \tag{4.7}$$

**Proposition 4.7.** *If* $\begin{pmatrix} r_1 \\ r_2 \end{pmatrix}$ *is an eigenvector of* $M$ (4.7) *corresponding to an eigenvalue* $\sigma$*, then it must satisfy* $r_1 = \sigma r_2$*. Moreover, we have*

(i) $r_2$ *is an eigenvector of* $G$ *associated to an eigenvalue* $\eta$*, where* $\eta$ *and* $\sigma$ *satisfy the relation*

$$\sigma^2 - \big((a-b) + (1+b)\eta\big)\sigma + (a-b) + b\eta = 0. \tag{4.8}$$

(ii) *Given any* $(a,b) \in [0,1[^2$*, then* $\rho(M) < 1$ *if, and only if,*

$$\frac{2(b-a)-1}{1+2b} < \underline{\eta}. \tag{4.9}$$

**Remark 4.8.** Though $G$ has $n$ eigenvalues, it can be shown that, given $a$ and $b$, $\rho(M)$ is determined only by $\underline{\eta}$ and $\overline{\eta}$. These extreme eigenvalues lie in $]-1,1[$ ($\gamma \in ]0, 2\beta[$) or even in $[0,1[$ ($\gamma \in ]0,\beta]$) by Lemma 4.4(i)(b)-(c).

**General partly smooth case** When $R$ is a general partly smooth function, then $U$ is nontrivial, and the spectral analysis of (4.3) becomes a generalized eigenvalue problem which is much more complex. Therefore, we assume $b = a$, in which case $M$ reads

$$M = \begin{bmatrix} (1+a)WG, & -aWG \\ \mathrm{Id}, & 0 \end{bmatrix}. \tag{4.10}$$

We have the following corollary of Proposition 4.7.

**Corollary 4.9.** *Let* $b = a$*. If* $\begin{pmatrix} r_1 \\ r_2 \end{pmatrix}$ *be an eigenvector of* $M$ *corresponding to an eigenvalue* $\sigma$*, then it must satisfy* $r_1 = \sigma r_2$*. Moreover* $r_2$ *is an eigenvector of* $G$ *related to eigenvalue* $\eta$*, where* $\eta$ *and* $\sigma$ *satisfy the relation*

$$\sigma^2 - (1+a)\eta\sigma + a\eta = 0, \tag{4.11}$$

*and* $\rho(M) < 1$ *if, and only if,*

$$\frac{-1}{1+2a} < \underline{\eta}. \tag{4.12}$$

**Remark 4.10.** Condition (4.12) holds naturally for $\gamma \in ]0,\beta]$, since by Lemma 4.4(ii), for such $\gamma$, $\underline{\eta} \geq 0$.

### 4.3 Local linear convergence of FB-type methods

Now we are able present the local linear convergence result of FB-type method, and start with the case where $R$ is locally polyhedral around $x^\star$.

**Theorem 4.11.** *Suppose* (**H.1**)-(**H.3**) *hold, and an FB-type method generates a sequence* $x_k \to x^\star \in \mathrm{Argmin}(\Phi)$ *such that* $R$ *is* locally polyhedral *around* $x^\star$*,* $F$ *is* $C^2$ *near* $x^\star$*, and conditions* (ND)*,* (RI) *are satisfied. If moreover* (4.4) *and* (4.9) *hold, then* $(x_k)_{k\in\mathbb{N}}$ *converges locally linearly to* $x^\star$*. More precisely, given any* $\rho \in [\rho(M), 1[$*, there exists* $K > 0$ *and a constant* $C > 0$*, such that for all* $k \geq K$*, there holds*

$$\|x_k - x^\star\| \leq C\rho^{k-K}\|x_K - x^\star\|.$$

**Proof.** Combining Proposition 4.5, Proposition 4.7 and [52, Section 2.1.2, Theorem 1], leads to the claimed result. □

**Remark 4.12.** $\rho(M)$ is the optimal rate. Indeed, when $a_k \equiv a, b_k \equiv b$ and $\gamma_k \equiv \gamma$, the $o(\cdot)$ term vanishes in (4.5) and thus, $\rho = \rho(M)$.

Let's turn to the case where $R$ is a general partly smooth function, but $b = a \in [0, \bar{a}]$ as in (4.10).

**Theorem 4.13.** *Suppose assumptions* (**H.1**)-(**H.3**) *hold, and the FB-type methods generate a sequence* $x_k \to x^\star \in \mathrm{Argmin}(\Phi)$ *such that* $R \in \mathrm{PSF}_{x^\star}(\mathcal{M}_{x^\star})$, $F$ *is* $C^2$ *near* $x^\star$, *and conditions* (ND), (RI) *are satisfied. If moreover* (4.4) *holds with* $b = a$, *and* (4.12) *is satisfied, then* $(x_k)_{k\in\mathbb{N}}$ *converges locally linearly to* $x^\star$. *More precisely, given any* $\rho \in [\rho(M), 1[$, *there exists* $K > 0$ *and a constant* $C > 0$, *such that for all* $k \geq K$, *there holds*

$$\|x_k - x^\star\| \leq C\rho^{k-K}\|x_K - x^\star\|.$$

**Proof.** This follows by combining Proposition 4.5, Corollary 4.9 and [52, Section 2.1.2, Theorem 1]. ☐

**Remark 4.14.**

(i) The limit $b = a$ in (4.4) does not mean that we should set $b_k = a_k, \forall k \in \mathbb{N}$ along the iterations.

(ii) In contrast to our previous work [39], which addresses the case of FB method, the rate estimates that we provide here are much sharper in general, and both estimates only coincide when $R$ is locally polyhedral (see the numerical experiments for more details). The main reasons underlying this is that, here, our rate estimate relies on the locally linearized iteration in Proposition 4.5 and the spectral properties of $M$, which takes intro account the geometry of the identified submanifold (its curvature for instance). This is not the case in our former work.

(iii) The obtained results can be readily extended to the variable metric FB splitting method [22], where a rate under an appropriate metric can be obtained. However for the sake of brevity, we do not pursue this further.

(iv) In our proof of local linear convergence, convexity does play a crucial role. For instance, it was only needed to show that the matrix $U$ is positive semi-definite. This suggests that our local linear convergence claims can be extended to the non-convex case, provided that the Riemannian Hessian of $R$ is assumed positive semi-definite at $x^\star$. In addition, to guarantee finite identification in the non-convex setting, we need global convergence of iFB to a critical point, which can be ensured if for instance $\Phi$ satisfies the (non-smooth) Kurdyka-Łojasiewicz inequality [15]. This will be left to a forthcoming paper.

The restricted injectivity condition (RI) plays an important role in our local convergence rate analysis and in general cannot be relaxed. However, for some special cases, such as when $R$ is locally polyhedral, it can be removed, at the price of less sharp rate estimation. This is formalized in the following statement.

**Theorem 4.15.** *Suppose that* (**H.1**)-(**H.3**) *hold, and an FB-type method creates a sequence* $x_k \to x^\star \in \mathrm{Argmin}(\Phi)$ *such that* $R$ *is locally polyhedral around* $x^\star$, $F$ *is* $C^2$ *near* $x^\star$, *and condition* (ND) *holds. If moreover there exists* $\epsilon > 0$ *and a subspace* $V$ *such that*

$$\ker\big(\mathrm{P}_{T_x}\nabla^2 F(x)\mathrm{P}_{T_x}\big) = V, \quad \forall x \in \mathbb{B}_\epsilon(x^\star) \cap (x^\star + T_{x^\star}).$$

*Then* $(x_k)_{k\in\mathbb{N}}$ *converges locally linearly to* $x^\star$.

The expressions of the local rate can be found by inspecting the proof.

## 4.4 Discussion

In this part, we present some discussions on the obtained local linear convergence result, and mainly focus on the difference FISTA and the iFB methods.

**FB is locally faster than FISTA** For the sake of brevity (the same conclusions hold true in the general case), we consider $b_k = a_k \equiv a \in [0,1]$ and $\gamma_k \equiv \gamma \in ]0, \beta]$ is fixed, in which case $\overline{\eta} \geq \underline{\eta} \geq 0$ (see Lemma 4.4(ii)), and thus condition (4.12) is in force. Moreover $\overline{\eta}$ is also the local convergence rate of the FB method, and $\rho(M)$ depends solely on $\overline{\eta}$ and the value of $a$. Recall that $\rho(M)$ is the best local linear convergence rate (see Theorem 4.13 and 4.11).

Figure 3 shows $\rho(M)$ as a function of $a$ for fixed $\overline{\eta}$. One can make the following observations:

(1) When $a \in [0, \overline{\eta}]$, we have $\rho(M) \leq \overline{\eta}$. This entails that if iFB is used with such a choice of inertial parameter, it will converges locally lineally faster than FB. For $a \in [\overline{\eta}, 1]$, the situation reverses as $\rho(M) \geq \overline{\eta}$, and iFB becomes slower than FB.

(2) In particular, as $a = 1$ for FISTA, we have $\rho(M) = \sqrt{\overline{\eta}} > \overline{\eta}$. In plain words, though FISTA is known to be globally faster (in terms of the objective) than FB, attaining the optimal $O(1/k^2)$ rate, locally, the situation radically changes as FISTA will always ends up being locally slower than FB. A similar observation is made in [58] for the special case of a variant of FISTA used to solve the LASSO problem. This explains in particular why many authors [27, 48] resort to restarting to accelerate local convergence of FISTA, which consists in resetting periodically the scheme to $a = 0$ which is more favorable to FISTA. Our predictions in Figure 3 gives clues on when to restart (*i.e.* detect the point in red on the rate curve). We will elaborate more on this in the numerical simulations in Section 5.5.

(3) $\rho(M)$ attains its minimal value at $a = \frac{(1-\sqrt{1-\overline{\eta}})^2}{\overline{\eta}}$, and this is the best convergence rate that can be achieved locally for FB-type methods.

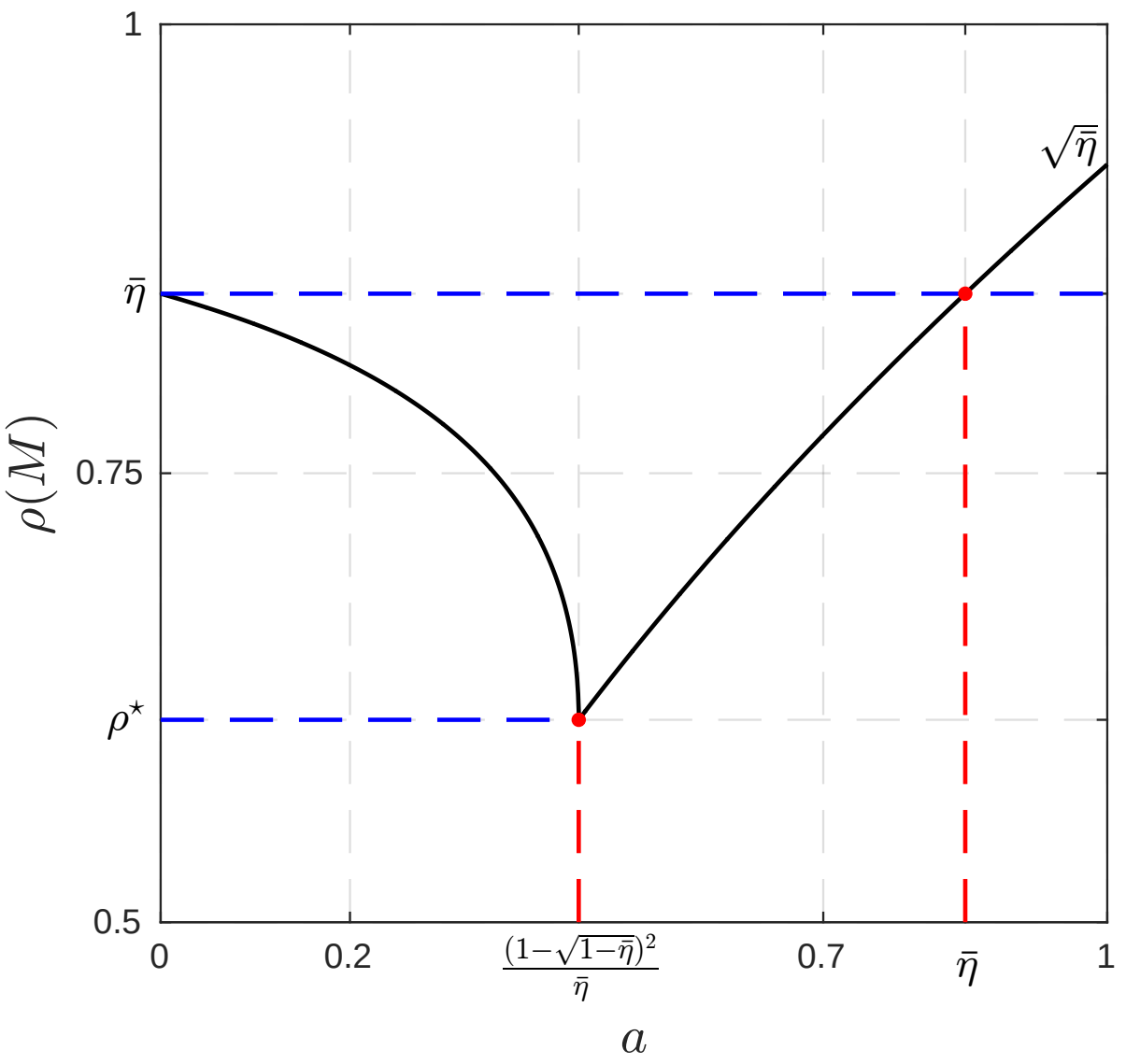

Figure 3: Let $b = a$, and assume $\underline{\eta}, \overline{\eta}$ are known and also close enough such that the spectral radius $\rho(M)$ is only affected by $\overline{\eta}$, then $\rho(M)$ is a function of $a$.

**Oscillation of the FISTA method** A typical feature of the FISTA method is that it is not monotone and locally oscillates [13], which makes the local convergence even slower, see Figure 4 or 5 for example, or [58] for a FISTA-variant applied to the LASSO problem. In fact, the iFB scheme shares this property as well when the inertial parameters are large.

Such oscillatory behaviour is due to the fact that, for those inertial parameters, the eigenvalue $\sigma_{\max}$ such that $|\sigma_{\max}| = \rho(M)$ is complex. It can then be seen that the oscillation period of $\|x_k - x^\star\|$ is exactly $\frac{\pi}{\theta}$, where $\theta$ is the argument of $\sigma_{\max}$. For the parameter settings used in Figure 3, *i.e.* $b = a$ and $\gamma \in ]0, \beta]$, we have

$$\begin{cases} a \in [0, \frac{(1-\sqrt{1-\overline{\eta}})^2}{\overline{\eta}}] : \sigma_{\max} \text{ is real,} \\ a \in ]\frac{(1-\sqrt{1-\overline{\eta}})^2}{\overline{\eta}}, 1] : \sigma_{\max} \text{ is complex,} \end{cases}$$

then as long as $a > \frac{(1-\sqrt{1-\overline{\eta}})^2}{\overline{\eta}}$, the iFB method locally oscillates. See Figure 5 for an example.

### 4.5 Acceleration

The finite time activity identification property (Theorem 3.4) implies that, the globally convex but non-smooth problem eventually becomes locally $C^2$-smooth, but possibly non-convex, constrained on the activity manifold. This opens the door to acceleration, and even finite termination, exploiting the structure of the objective and that of the identified manifold. There are several ways to achieve this goal as we explain hereafter.

**Optimal first-order method** In this case, the idea is to keep the scheme implemented in Algorithm 1, and to refine the parameters to minimize the local convergence rate established in Section 4. Indeed, as shown in Figure 3 and the discussion that follows, there is a proper choice of the inertial parameters $a$ and $b$ that minimizes $\rho(M)$. More precisely, choose $\gamma \in ]0, \beta]$, then $\overline{\eta} = 1 - \alpha\gamma \geq \underline{\eta} \geq 1 - \gamma/\beta \geq 0$, and $\rho(M)$ depends only on $\overline{\eta}$, $a$ and $b$. Then with fixed $\gamma$ (hence $\overline{\eta}$), $\rho(M)$ attains its minimal value for $a$ and $b$ satisfying

$$\begin{cases} b = a : a = \dfrac{(1-\sqrt{1-\overline{\eta}})^2}{\overline{\eta}} = \dfrac{1-\sqrt{\alpha\gamma}}{1+\sqrt{\alpha\gamma}}, \\ b \neq a : a = \big(1-\sqrt{1-\overline{\eta}}\big)^2 + b(1-\overline{\eta}) = (1-\sqrt{\alpha\gamma})^2 + b\alpha\gamma, \end{cases} \tag{4.13}$$

and the optimal value $\rho^\star$ of $\rho(M)$ reads

$$\rho^\star = 1 - \sqrt{1-\overline{\eta}} = 1 - \sqrt{\gamma\alpha}, \tag{4.14}$$

where the second equality comes from (4.2) and Lemma 4.4. This is a decreasing function of $\gamma$, and $\rho^\star = 1 - \sqrt{\alpha\beta}$ is then the minimal rate attained for $\gamma = \beta$. This rate is in agreement with that [46, Theorem 2.2.2]. If one can afford $\gamma \geq \beta$ as in our iFB schemes, owing to the result of [52, Section 3.2.1], the best local linear rate is actually

$$\underline{\rho}^\star = \frac{1-\sqrt{\alpha\beta}}{1+\sqrt{\alpha\beta}} \quad \text{for} \quad \gamma = \frac{4\beta}{(1+\sqrt{\alpha\beta})^2}, \quad a = \Big(\frac{1-\sqrt{\alpha\beta}}{1+\sqrt{\alpha\beta}}\Big)^2 \quad \text{and} \quad b = 0.$$

This is known to be the optimal rate that matches the lower complexity bounds for first-order methods to solve the class of problems ($\mathcal{P}_{\mathrm{opt}}$) if $F$ were also $\alpha$-strongly convex [46, Theorem 2.1.13]. In comparison, for the FB method (*i.e.* $a = b = 0$), the optimal rate is $\rho^\star = \overline{\eta}^\star = \frac{1-\alpha\beta}{1+\alpha\beta}$ attained for $\gamma = \frac{2\beta}{1+\alpha\beta}$.

**Finite convergence in the polyhedral case** Finite termination can be obtained if $R$ is locally polyhedral around $x^\star$, and $F$ is quadratic, *i.e.* problem ($\mathcal{P}_\lambda$) with $R$ locally polyhedral around $x^\star$. In this situation, under hypothesis (ND), we have finite identification of $x^\star + T_{x^\star}$. In addition, (RI) is equivalent to injectivity of the linear operator $L$ on $T_{x^\star}$. Altogether, this allows to show that $x^\star$ can be written explicitly as

$$x^\star = L^{*,+}_{T_{x_K}} y - \lambda\big(L^*_{T_{x_K}} L_{T_{x_K}}\big)^+ \mathrm{P}_{T_{x_K}}\big(\partial R(x_K)\big) - L\mathrm{P}_{x_K + T_{x_K}}(0),$$

for $K$ sufficiently large. When the manifold is linear, *i.e.* $x^\star \in T_{x^\star}$, the last term vanishes and the above relation can be implemented in practice.

**High-order acceleration: Newton method** Once the activity manifold has been identified, one can switch to Newton-type methods for locally minimizing $\Phi$. This can be done either using local parameterizations obtained from $\mathcal{U}$-Lagrangian theory or from Riemannian geometry [36, 42, 56]. One can also use the Riemannian version of the non-linear conjugate gradient method [56]. For these schemes, one can also show respectively quadratic and superlinear convergence since $\nabla^2_{\mathcal{M}_{x^\star}} \Phi(x^\star)$ is positive definite by Proposition 4.1(ii).

# 5 Numerical experiments

In this section, we illustrate the obtained results by some popular examples drawn from linear inverse problems in signal processing and machine learning (including sparse recovery). We first start by discussing a few examples of partly smooth functions that are widely used in those applications.

## 5.1 Examples of partly smooth functions

**Example 5.1 ($\ell_1$-norm).** For $x \in \mathbb{R}^n$, the $\ell_1$-norm is defined as

$$R(x) = \|x\|_1 \stackrel{\text{def}}{=} \textstyle\sum_{i=1}^n |x_i|,$$

which is polyhedral, hence partly smooth at any $x$ relative to the subspace

$$\mathcal{M} = T_x \stackrel{\text{def}}{=} \big\{u \in \mathbb{R}^n : \mathrm{supp}(u) \subseteq \mathrm{supp}(x)\big\}, \ \mathrm{supp}(x) \stackrel{\text{def}}{=} \{i : x_i \neq 0\}.$$

Its Riemannian gradient at $x$ is $\mathrm{sign}(x_i)$ for $i \in \mathrm{supp}(x)$, and 0 otherwise. Its Riemannian Hessian vanishes.

**Example 5.2 ($\ell_{1,2}$-norm).** Let the index set $\{1, \ldots, n\}$ be partitioned into non-overlapping blocks $\mathcal{B}$ such that $\bigcup_{b\in\mathcal{B}} b = \{1, \ldots, n\}$. The $\ell_{1,2}$-norm of $x$ is given by

$$R(x) = \|x\|_{1,2} \stackrel{\text{def}}{=} \textstyle\sum_{b\in\mathcal{B}} \|x_b\|,$$

where $x_b = (x_i)_{i\in b} \in \mathbb{R}^{|b|}$. Though this function is not polyhedral, it is easy to see that it is partly smooth at $x$ relative to the subspace

$$\mathcal{M} = T_x \stackrel{\text{def}}{=} \big\{u \in \mathbb{R}^n : \mathrm{supp}_{\mathcal{B}}(u) \subseteq \mathcal{S}_{\mathcal{B}}\big\}, \ \mathcal{S}_{\mathcal{B}} \stackrel{\text{def}}{=} \bigcup \{b : x_b \neq 0\}.$$

It is straightforward to show that

$$\big(\nabla_{\mathcal{M}} \|x\|_{1,2}\big)_b = \begin{cases} x_b/\|x_b\| & \text{if } x_b \neq 0 \\ 0 & \text{otherwise} \end{cases} \ \text{and} \ \nabla^2_{\mathcal{M}} \|x\|(x) = \delta_x \circ Q_{x^\perp},$$

where,

$$\delta_x : T_x \to T_x, v \mapsto \begin{cases} v_b/\|x_b\| & \text{if } x_b \neq 0, \\ 0 & \text{otherwise,} \end{cases} \ \text{and} \ \ Q_{x^\perp} : T_x \to T_x, v \mapsto \begin{cases} v_b - \frac{\langle x_b, v_b\rangle}{\|x_b\|^2} x_b & \text{if } x_b \neq 0, \\ 0 & \text{otherwise.} \end{cases}$$

**Example 5.3 (Total Variation).** If $R_0 \in \mathrm{PSF}_{D^*x}(\mathcal{M}_0)$, then, under a mild transversality condition, it is shown in [37, Theorem 4.2] that $R \in \mathrm{PSF}_x(\mathcal{M})$ where $\mathcal{M} = \{u \in \mathbb{R}^n : D^*u \in \mathcal{M}_0\}$. Popular examples include the anisotropic total variation (TV) semi-norm in which case $R_0 = \|\cdot\|_1$ and $D^* = D_{\mathrm{DIF}}$ is a finite difference approximation of the derivative [55]. For TV, $R$ is then polyhedral, hence partly smooth at $x$ relative to

$$\mathcal{M} = T_x \stackrel{\text{def}}{=} \{u \in \mathbb{R}^n : \mathrm{supp}(D^*u) \subseteq \mathrm{supp}(D^*x)\}.$$

Its Riemannian gradient reads $\mathrm{P}_{T_x}\mathrm{sign}(D^*x)$ and its Riemannian Hessian vanishes.

**Example 5.4 ($\ell_\infty$-norm).** For $x \in \mathbb{R}^n$, the anti-sparsity promoting $\ell_\infty$-norm is defined as following

$$R(x) = \|x\|_\infty \stackrel{\text{def}}{=} \max_{1 \le i \le n} |x_i|.$$

It can verified that $R$ is a polyhedral norm, hence partly smooth at $x$ relative to

$$\mathcal{M} = T_x \stackrel{\text{def}}{=} \mathbb{R}s_{I(x)},\ \ I(x) \stackrel{\text{def}}{=} \{i :\ |x_i| = \|x\|_\infty\}, s_i \stackrel{\text{def}}{=} \begin{cases} \mathrm{sign}(x_i), & \text{if } i \in I(x), \\ 0, & \text{otherwise.} \end{cases}$$

The Riemannian gradient of $\|\cdot\|_\infty$ at $x$ is $s/|I(x)|$, and its Riemannian Hessian vanishes.

**Example 5.5 (Nuclear norm).** For $x \in \mathbb{R}^{n_1 \times n_2}$ with $\mathrm{rank}(x) = r$, let $x = U\mathrm{diag}(\sigma(x))V^*$ be a reduced rank-$r$ SVD decomposition, where $U \in \mathbb{R}^{n_1 \times r}$ and $V \in \mathbb{R}^{n_2 \times r}$ have orthonormal columns, and $\sigma(x) \in (\mathbb{R}_+ \setminus \{0\})^r$ is the vector of singular values $(\sigma_1(x), \cdots, \sigma_r(x))$ in non-increasing order. Low-rank is the spectral extension of vector sparsity to matrix-valued data $x \in \mathbb{R}^{n_1 \times n_2}$, *i.e.* imposing sparsity on the singular values of $x$. The nuclear norm is thus defined as

$$R(x) = \|x\|_* \stackrel{\text{def}}{=} \|\sigma_i(x)\|_1.$$

Piecing together [23, Theorem 3.19] and Example 5.1, the nuclear norm can be shown to be partly smooth at $x$ relative to the set of fixed-rank matrices

$$\mathcal{M} \stackrel{\text{def}}{=} \{z \in \mathbb{R}^{n_1 \times n_2} : \mathrm{rank}(z) = r\},$$

which is a $C^2$-manifold around $x$ of dimension $(n_1 + n_2 - r)r$, see [35, Example 8.14].

Moreover, we have

$$T_x = \{UA^* + BV^* :\ A \in \mathbb{R}^{n_2 \times r}, B \in \mathbb{R}^{n_1 \times r}\} \text{ and } \nabla_{\mathcal{M}}\|x\|_* = UV^*.$$

From [60, Example 21], one can show that for $h \in T_x$,

$$\nabla^2_{\mathcal{M}}\|x\|_*(h) = \mathrm{P}_{T_x}\nabla^2\widetilde{\|x\|_*}(\mathrm{P}_{T_x}h),$$

where

$$\widetilde{\|z\|_*} = \widetilde{\|\sigma(z)\|_1} = \textstyle\sum_{i=1}^r \sigma_i(z),$$

is a $C^2$-smooth (and even convex) representation of the nuclear norm on $\mathcal{M}$ near $x$, obtained owing to the smooth transfer principle [23, Corollary 2.3]. The expression of the (Euclidian) Hessian $\nabla^2\widetilde{\|z\|_*}$ can be obtained in several ways, see [60, Example 21] for details.

### 5.2 Linear inverse problems

In this part, we apply our results to the setting of linear inverse problems. Consider the following forward observation of a vector $x_{\rm ob} \in \mathbb{R}^n$

$$y = Lx_{\rm ob} + w, \tag{5.1}$$

where $y \in \mathbb{R}^m$ is the observation, $L : \mathbb{R}^n \to \mathbb{R}^m$ is some linear operator, and $w \in \mathbb{R}^m$ stands for noise. Solving such linear inverse problems can be cast as the optimization problem

$$\min_{x\in\mathbb{R}^n} \frac{1}{2}\|y - Lx\|^2 + \lambda R(x), \tag{$\mathcal{P}_\lambda$}$$

where $\lambda > 0$ is the regularization parameter, $R \in \Gamma_0(\mathbb{R}^n)$ encodes prior knowledge on $x_{\rm ob}$ and hence promotes objects similar to to it, and $\lambda > 0$ is a regularization parameter. Moreover, when there is no noise in the observation (5.1), namely $w = 0$, the following equality constrained problem should be considered

$$\min_{x\in\mathbb{R}^n} R(x) \ \text{ s.t. } \ Lx = Lx_{\rm ob}. \tag{$\mathcal{P}_0$}$$

The following result is a straightforward generalization of [62, Theorem 1] to any FB-type method, using Theorem 3.4 and Theorem 4.11 (or 4.13).

**Proposition 5.6.** *Assume that $R \in \mathrm{PSF}_{x_{\rm ob}}(\mathcal{M}_{x_{\rm ob}})$, and condition*

$$\ker(L) \cap T_{x_{\rm ob}} = \{0\} \ \textit{and} \ \big(L^+_{T_{x_{\rm ob}}} L\big)^T \nabla_{\mathcal{M}_{x_{\rm ob}}} R(x_{\rm ob}) \in \mathrm{ri}\big(\partial R(x_{\rm ob})\big), \tag{5.2}$$

*hold. If moreover $w$ is sufficiently small and $\lambda$ is chosen in the order of $\|w\|$, then ($\mathcal{P}_\lambda$) admits a unique solution $x^\star$ with $\mathcal{M}_{x^\star} = \mathcal{M}_{x_{\rm ob}}$, and the FB-type methods will identify $\mathcal{M}_{x^\star}$ in finite time, and then converge locally linearly.*

This proposition implies that under the given conditions, the minimizer of ($\mathcal{P}_\lambda$) lies in the same manifold as the feasible point of ($\mathcal{P}_0$). It is now sufficient to infer when (5.2) is satisfied for the above proposition to hold true. For instance, when $L$ is a random Gaussian measurement matrix, nice and easily verifiable conditions can be stated for the examples introduced in Section 5.1 above.

**Proposition 5.7.** *Choose $L$ from the standard Gaussian ensemble, i.e. the entries of $L$ are independent copies of a mean-zero and standard Gaussian random variable. Then (5.2) is in force with high probability in the following cases:*

(i) $R = \|\cdot\|_1$*: let $s = \|x_{\rm ob}\|_0$, if $m \geq 2cs\log(n) + s$ for some $c > 1$;*

(ii) $R = \|\cdot\|_{1,2}$*: let $s$ be the number of non-zero blocks, if $m \geq (1+c)s\big(\sqrt{n/N_{\mathcal{B}}} + \sqrt{2\log(N_{\mathcal{B}})}\big)^2 + sn/N_{\mathcal{B}}$ where $c > 1$, and $N_{\mathcal{B}}$ is the total number of blocks;*

(iii) $R = \|\cdot\|_\infty$*: let $I(x) = \{i : \ |(x_{\rm ob})_i| = \|x_{\rm ob}\|_\infty\}$ and $s = |I(x)|$, if $m \geq n - s + 2cs\log(s/2)$, where $c > 1$;*

(iv) $R = \|\cdot\|_*$*: let $r = \mathrm{rank}(x_{\rm ob})$, $x_{\rm ob} \in \mathbb{R}^{n_1\times n_2}$, if $m \geq cr(3n_1 + 3n_2 - 5r)$ for some $c > 1$.*

**Proof.** This follows from [18, Section 3] for (i), (ii) and (iv), and (iii) from [61, Theorem 7]. □

### 5.3 Experiments setup

**Recovery from random measurements** We consider solving ($\mathcal{P}_\lambda$) with $R$ being $\ell_1, \ell_{1,2}, \ell_\infty$-norms, TV semi-norm and nuclear norm. The observations are generated according to (5.1). Here $L$ is generated from the standard Gaussian ensemble and the following parameters:

| | |
|---|---|
| **$\ell_1$-norm** | $(m, n) = (48, 128)$, $\|x_{\mathrm{ob}}\|_0 = 8$; |
| **$\ell_{1,2}$-norm** | $(m, n) = (60, 128)$, $x_{\mathrm{ob}}$ has 3 non-zero blocks of size 4; |
| **$\ell_\infty$-norm** | $(m, n) = (123, 128)$, $\lvert I(x_{\mathrm{ob}})\rvert = 10$; |
| **Total Variation** | $(m, n) = (48, 128)$, $\|D_{\mathrm{DIF}} x_{\mathrm{ob}}\|_0 = 8$ where $D_{\mathrm{DIF}}$ is the finite difference operator; |
| **Nuclear norm** | $(m, n) = (1425, 2500)$, $x_{\mathrm{ob}} \in \mathbb{R}^{50\times 50}$ and $\mathrm{rank}(x_{\mathrm{ob}}) = 5$. |

It can be noticed that the number of measurements $m$ is chosen sufficiently large such that Proposition 5.7 allows to assert that (ND) and (RI) are verified at $x_{\mathrm{ob}}$. We also choose $\|w\|$ small enough and $\lambda$ in the order of $\|w\|$ so that Proposition 5.6 applies.

**TV deconvolution** We also consider a 2D image processing problem, where $y$ is a degraded image generated according to (5.1), where $L$ is a circular convolution matrix with a Gaussian kernel. The (anisotropic) TV regularizer (see Example 5.3), which is polyhedral, is used.

Note however that for a sparse deconvolution problem through $\ell_1$-minimization, Proposition 5.7 does not apply, hence entailing that exact recovery of the support of $x_{\mathrm{ob}}$ in general is impossible, see [25]. However, under the same conditions on $x_{\mathrm{ob}}$ and $\lambda$ as in Proposition 5.7, $x^\star$ has a support slightly larger than that of $x_{\mathrm{ob}}$, and moreover, $x^\star$ satisfies both (ND) and (RI). See [25, Corollary 1].

### 5.4 Comparison of the FB-type methods

**Parameter settings** For all the methods in comparison (FB, iFB and FISTA), we fix $\gamma_k \equiv \beta$. For the sequence convergent FISTA method [19, 9], two different choices of $q$ are considered, which are 2 and 50. For the iFB method, we let $b_k = a_k$, and use the following rule to update $a_k$. Let $t_0 = 1$, $p \in ]0, +\infty[$, then

$$t_k = \frac{1 + \sqrt{1 + p t_{k-1}^2}}{2}, \quad a_k = \frac{t_{k-1} - 1}{t_k} \begin{cases} p \in ]0, 4[ : t_k \to \dfrac{4}{4-p}, \ a_k \to \dfrac{p}{4}, \\ p \in [4, +\infty[ : t_k \to +\infty, \ a_k \to \dfrac{2}{\sqrt{p}}. \end{cases} \tag{5.3}$$

In this test we choose $p = 4(\sqrt{5} - 2 - 10^{-3})$ so that Theorem 2.3 applies. Note that in the original FISTA paper [14], (5.3) is also used but with $p = 4$ fixed.

The convergence profiles of $\|x_k - x^\star\|$ are shown in Figure 4. As demonstrated by all the plots, identification and local linear convergence occurs after finite time. The solid lines (denoted as "P") represent the observed profiles, while dashed ones (denoted as "T") stand for the theoretically predicted ones. The positions of the cyan points (or the starting points of the dashed lines) stand for the iteration at which $\mathcal{M}_{x^\star}$ has been identified.

**Tightness of predicted rates** For the $\ell_1, \ell_\infty$-norms and TV semi-norm, our predicted rates coincide exactly with the observed ones (same slopes for the dashed and solid lines). This is due to the fact that they are all polyhedral and $F$ is quadratic. Note that for FISTA, which is non-monotone, the prediction coincides with the envelope of the oscillations. For the $\ell_{1,2}$-norm, though it is not polyhedral, our predicted rates still are very tight, due to the fact that the Riemannian Hessian is taken into account. Then for the nuclear norm, whose active manifold is not anymore a subspace, our estimation becomes slightly less sharp compared to the other examples, though barely visible on the plots. For both the $\ell_{1,2}$-norm and nuclear norm, since the Riemannian Hessian is taken into account, the predicted rates are much sharper than our previous estimates for the FB method in [39].

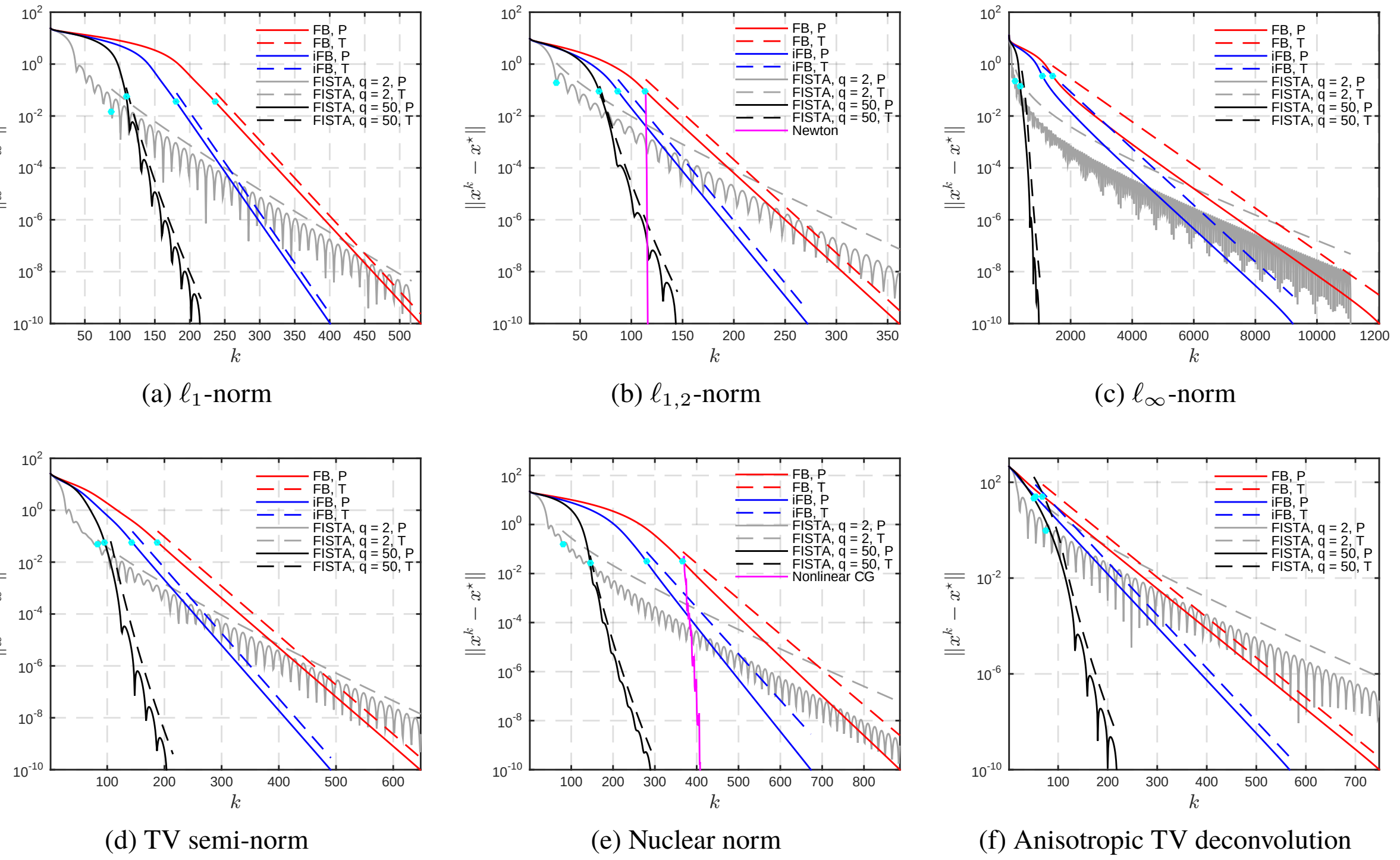

(a) $\ell_1$-norm (b) $\ell_{1,2}$-norm (c) $\ell_\infty$-norm

(d) TV semi-norm (e) Nuclear norm (f) Anisotropic TV deconvolution

Figure 4: Local linear convergence and comparison of the FB-type methods (FB, iFB and FISTA) in terms of $\|x_k - x^\star\|$. We fix $\gamma_k \equiv \beta$ for all the methods, moreover, for the iFB method, we let $b_k = a_k \equiv \sqrt{5} - 2 - 10^{-3}$, and for the FISTA method, $q = 2, 50$ are considered. For each figure, "P" stands for practical observed profiles, while "T" indicates theoretical predictions. The cyan points indicate the iteration at which $\mathcal{M}_{x^\star}$ has been identified.

For the image deconvolution problem, assumptions (ND) and (RI) are checked a posteriori (verified for this experiment). This together with the fact that the anisotropic TV is polyhedral justifies that the predicted rate is again exact (up to machine precision).

**Comparison of the methods** From the numerical results, we can draw the following remarks:

(i) Overall, FISTA with $q = 50$ (black line) is the fastest while $q = 2$ (gray line) is the slowest. FB and iFB are sandwiched between them with iFB being the faster one.

(ii) For the finite activity identification, however, FISTA $q = 2$ in general shows the fastest identification (see the starting points of the dashed lines), and FB is the slowest.

(iii) Locally, similar to the global convergence, FISTA $q = 50$ has the fastest rate and $q = 2$ is the slowest. Again, FB and iFB are between them with iFB being faster than FB.

It can be concluded from the above remarks that, in practice, FISTA method with $q = 2$ is not a wise choice if high accuracy solutions are needed. Indeed, under this choice, $a_k$ converges to 1 too fast, and this hampers its local behaviour as the discussions we anticipated in Section 4.4 (see Figure 3). In fact, such behaviour of $a_k$ can be avoided by choosing relatively bigger $q$, and this is exactly what the difference between $q = 2$ and $q = 50$ implies. In our tests, $q \in [50, 100]$ seems to a good trade-off, even bigger $q$ is not recommended since it may lead to a much slower activity identification. A similar observation is also mentioned in [19], where the authors only tried $q = 2, 3, 4$. It should be noted that the original FISTA method [14] has almost the same behaviour as the case $q = 2$.

It should be pointed out that the local rate of FISTA $q = 50$ being faster than FB does not contradict with our claim in Section 4.4 that FB is faster than FISTA locally. The reason is that we are limited by machine accuracy, and bigger value of $q$ delays the speed at which $a_k$ approaches to 1 which actually makes FISTA behaviour similar to the iFB method.

**High-order acceleration** For the $\ell_1, \ell_\infty$-norms and TV semi-norm, since they are polyhedral, finite termination can be obtained once the manifold is identified. For $\ell_{1,2}$-norm which is not polyhedral, we applied the Riemannian Newton method which converges quadratically, leading to a dramatic acceleration as can be seen in Figure 4(b). For the nuclear norm, a non-linear conjugate gradient method is applied, leading again to a much faster (super-linear) local convergence.

**Oscillation of the FISTA method** As observed from Figure 4, FISTA method oscillates for both choices of $q$. No oscillation appears for the iFB method since the value of the inertial parameter is not big enough. In order to have a better visualisation of the oscillation of iFB/FISTA methods, we choose the LASSO problem for illustration, set $b = a$ and locally adjust the value of $a$ so that the oscillation period is integer. The result is shown in Figure 5, where the oscillation period of the tested example is 20.

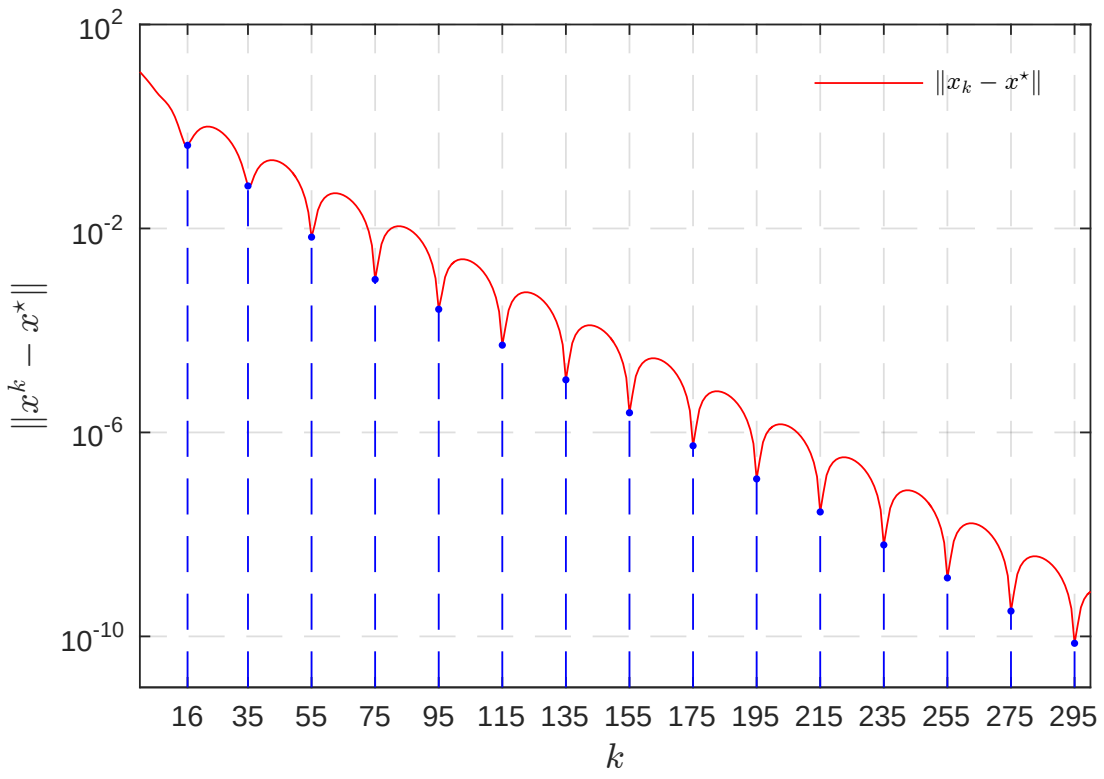

Figure 5: Local oscillation of the iFB/FISTA methods on LASSO problem. Local oscillation of the iFB method, where the oscillation period is 20.

### 5.5 Comparisons of step-size and inertial settings

In this section, we provide more comparisons of the iFB method, on the choices of different step-size $\gamma_k$ and also the difference between the inertial parameters $a_k, b_k$ respectively.

**Comparison of $\gamma_k \equiv \beta$ vs $\gamma_k \equiv 1.5\beta$** We compare the difference between different step-sizes, and two choices of $\gamma_k$ are considered: $\gamma_k \equiv \beta$ and $\gamma_k \equiv 1.5\beta$, and the corresponding inertial parameter are,

- iFB $\gamma_k \equiv \beta$: $a_k = b_k \equiv \sqrt{5} - 2 - 10^{-3}$ same as above tests such that Theorem 2.3 applies;
- iFB $\gamma_k \equiv 1.5\beta$: $a_k, b_k$ are chosen according to (2.3) such that Theorem 2.1 applies.

In comparison, FISTA method with $q = 2$ and $q = 50$ is also added. Two numerical experiments on $\ell_1$-norm and nuclear norm are illustrated in Figure 6. From the numerical results, we can infer the following observations.

- For FB, larger $\gamma_k$ leads to faster global convergence and activity identification. However this does not mean that the bigger the better locally. As we discussed in Section 4.5, the best choice to get the optimal local linear rate is $2\beta/(1 + \alpha\beta)$.
- iFB is faster than FB under the same choice of $\gamma_k$. FISTA $q = 50$ is no longer the fastest one, while it is outperformed by iFB $\gamma_k \equiv 1.5\beta$ on the LASSO problem. Moreover, it should be noted that, the inertial parameters of iFB can be optimized according to (4.13), which can make the iteration even faster.

In accompany with the high-order acceleration result present above, it can be conclude that in practice, the *inertial+high-order method* hybrid strategy is an ideal choice for solving ($\mathcal{P}_{\text{opt}}$).

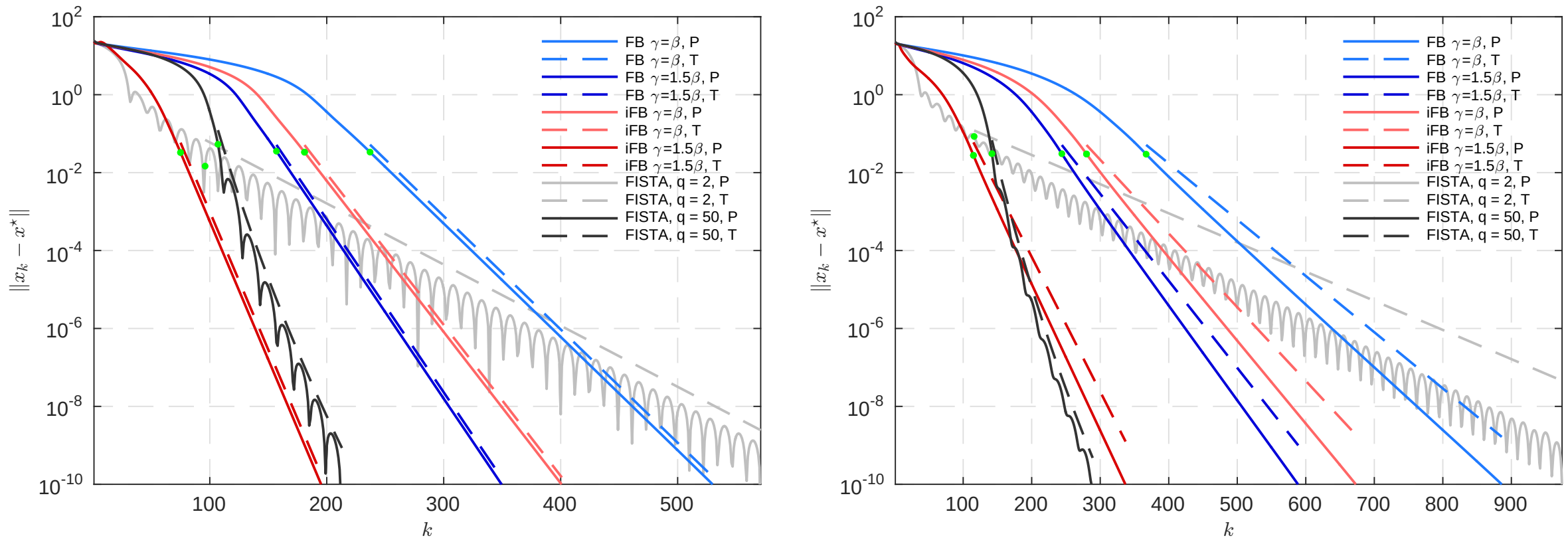

(a) $\ell_1$-norm

(b) Nuclear norm

Figure 6: Comparison of iFB method with different step-sizes.

**Comparison of $a_k$ vs $b_k$** Now let's assess the influence of inertial parameter choice, same as above $\ell_1$-norm and nuclear norm are considered. The step-size $\gamma$ is fixed as $\gamma_k \equiv \beta$.

For the iFB method, the online updated rule (2.3) is applied, with $c_{a,k} = c_{b,k} = \frac{10^5}{k^2\|x_k - x_{k-1}\|^2}$, and 4 different combinations of $(a, b)$ are considered, which are

$$(0.3, 0.2 \,\text{or}\, 0.6) \ \text{ and } \ (0.8, 0.2 \,\text{or}\, 0.6).$$

For both examples, if we let $b_k \equiv 0$, then the optimal local choice $a_{\text{opt}}$ obtained through (4.13) is between 0.3 and 0.8. The obtained plots are depicted in Figure 7, whence we summarize the following observations:

(i) The time to activity identification is more dependent on the value of $a$. Clearly, relatively bigger values of $a$ lead to a faster identification. On the other hand, when $a < a_{\text{opt}}$ (case $a = 0.3$), bigger values of $b$ lead to slower identification, while the opposite situation occurs when $a > a_{\text{opt}}$ (case $a = 0.8$).

(ii) The convergence rate also depends more on the choice of $a$, since with fixed $a$, the rate difference caused by different values of $b$ is small, see the blue dashed/solid lines, and magenta ones.

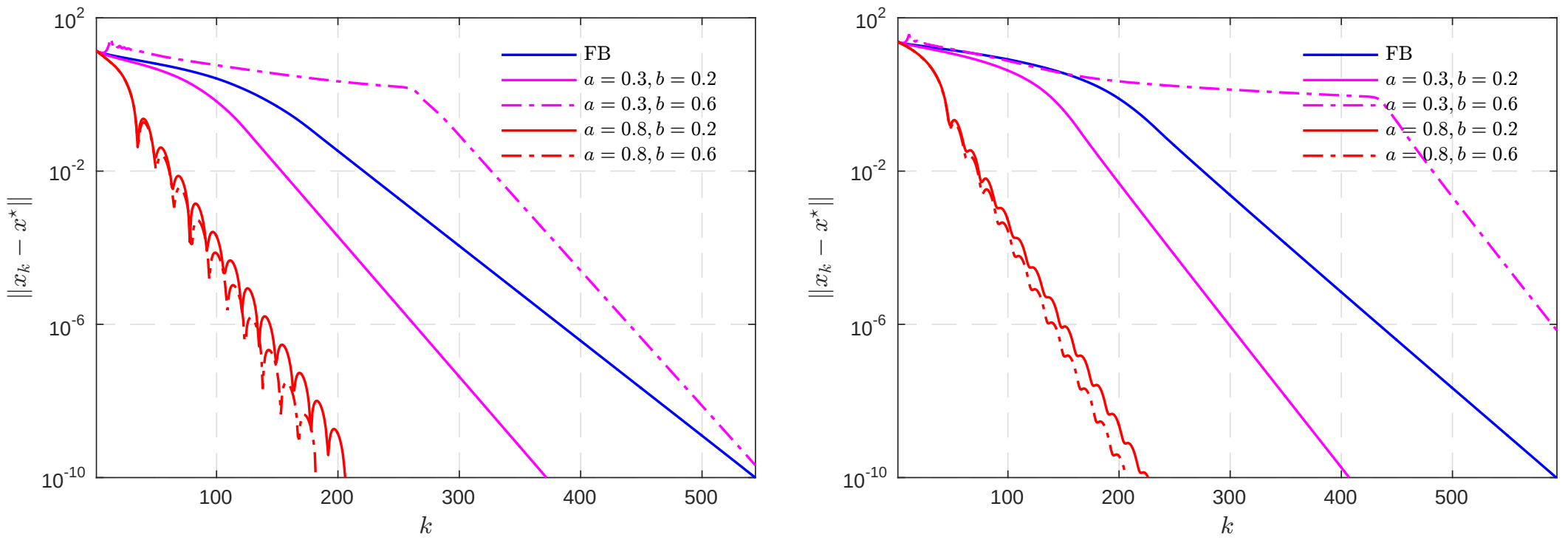

(a) $\ell_1$-norm (b) Nuclear norm

Figure 7: Comparisons on the difference between the inertial parameters $a_k$ and $b_k$, the step-size $\gamma$ is fixed as $\beta$.

# 6 Discussion and conclusion

In this paper, we proposed a generalized inertial Forward–Backward splitting scheme which covers several existing methods as special cases, and presented the corresponding global convergence analysis. Under partial smoothness, we established that this class of schemes identify the active manifold in finite time, and then converge locally linearly. The predicted rates were shown to be very sharp. We verified our theoretical findings with concrete numerical examples from signal/image processing and machine learning.

Most of our results can be extended to the non-convex setting by introducing appropriate supplementary assumptions, such as prox-regularity and the nonsmooth Kurdyka-Łojasiewicz inequality. This will be treated in a future work.

# Acknowledgements

This work has been partly supported by the European Research Council (ERC project SIGMA-Vision). JF was partly supported by Institut Universitaire de France.

# A Proofs of Section 2

Throughout this section, $\mathcal{H}$ denotes a real Hilbert space. We give a proof in the most general setting, *i.e.* solving ($\mathcal{P}_{\mathrm{inc}}$) on $\mathcal{H}$. We denote $\to$ strong convergence and $\rightharpoonup$ weak convergence on $\mathcal{H}$. We first briefly introduce some preliminaries which are needed for the convergence proof. Let $A : \mathcal{H} \rightrightarrows \mathcal{H}$ be a set-valued operator. The graph of $A$ is the set $\mathrm{gph}\, A = \{(x, y) \in \mathcal{H} \times \mathcal{H} | y \in A(x)\}$, and its zeros set is $\mathrm{zer} A = \{x \in \mathcal{H} | 0 \in A(x)\}$.

A set-valued operator $A : \mathcal{H} \rightrightarrows \mathcal{H}$ is monotone if

$$\big(\forall (x, v) \in \mathrm{gph}\, A\big), \big(\forall (y, u) \in \mathrm{gph}\, A\big), \ \langle x - y,\, v - u\rangle \geq 0. \tag{A.1}$$

It is moreover maximal monotone if $\mathrm{gph}\, A$ can not be contained in the graph of any other monotone operator.

Let $\beta \in ]0, +\infty[$, $B : \mathcal{H} \to \mathcal{H}$, then $B$ is $\beta$-cocoercive if the following holds

$$(\forall x, y \in \mathcal{H}),\ \beta\|Bx - By\|^2 \le \langle Bx - By,\, x - y\rangle, \tag{A.2}$$

which indicates that $B$ is $\beta^{-1}$-Lipschitz continuous.

**Proof of Theorem 2.1.** Let $x^\star \in \mathrm{zer}(A + B)$, *i.e.* a solution ($\mathcal{P}_{\text{inc}}$), which exists thanks to (**H.6**). From (1.4), we get

$$\begin{aligned} -B(x^\star) &\in A(x^\star),\\ (y_{a,k} - x_{k+1}) - \gamma_k B(y_{b,k}) &\in \gamma_k A(x_{k+1}). \end{aligned} \tag{A.3}$$

Define the following quantities

$$\varphi_k = \frac{1}{2}\|x_k - x^\star\|^2,\ E_{x,k} = \frac{1}{2}\|x_k - x_{k-1}\|^2,\ E_{a,k+1} = \frac{1}{2}\|y_{a,k} - x_{k+1}\|^2,\ E_{b,k+1} = \frac{1}{2}\|y_{b,k} - x_{k+1}\|^2. \tag{A.4}$$

By definition of $y_{a,k}$ we have

$$\begin{aligned} \varphi_k - \varphi_{k+1} &= \frac{1}{2}\langle x_k - x^\star,\, x_k - x^\star\rangle - \frac{1}{2}\langle x_{k+1} - x^\star,\, x_{k+1} - x^\star\rangle\\ &= \frac{1}{2}\langle x_k - x_{k+1} - x^\star + 2x_{k+1} - x^\star,\, x_k - x_{k+1}\rangle\\ &= E_{x,k+1} + \langle x_k - y_{a,k} + y_{a,k} - x_{k+1},\, x_{k+1} - x^\star\rangle\\ &= E_{x,k+1} + \langle y_{a,k} - x_{k+1},\, x_{k+1} - x^\star\rangle - a_k\langle x_k - x_{k-1},\, x_{k+1} - x^\star\rangle. \end{aligned} \tag{A.5}$$

Meanwhile, by virtue of the monotonicity of $A$ and (A.3), we have

$$\begin{aligned} \langle \gamma_k u_{k+1} - \gamma_k u^\star,\, x_{k+1} - x^\star\rangle &\ge 0, \quad \forall u_{k+1} \in A(x_{k+1}), u^\star \in A(x^\star)\\ \langle (y_{a,k} - x_{k+1}) - \gamma_k B(y_{b,k}) + \gamma_k B(x^\star),\, x_{k+1} - x^\star\rangle &\ge 0, \end{aligned}$$

which leads to

$$\langle y_{a,k} - x_{k+1},\, x_{k+1} - x^\star\rangle \ge \gamma_k\langle B(y_{b,k}) - B(x^\star),\, x_{k+1} - x^\star\rangle.$$

Combining this with (A.5), we obtain

$$\varphi_k - \varphi_{k+1} \ge E_{x,k+1} + \gamma_k\langle B(y_{b,k}) - B(x^\star),\, x_{k+1} - x^\star\rangle - a_k\langle x_k - x_{k-1},\, x_{k+1} - x^\star\rangle. \tag{A.6}$$

For $\langle x_k - x_{k-1},\, x_{k+1} - x^\star\rangle$, we have

$$\begin{aligned} \langle x_k - x_{k-1},\, x_{k+1} - x^\star\rangle &= \langle x_k - x_{k-1},\, x_{k+1} - x_k + x_k - x^\star\rangle\\ &= \langle x_k - x_{k-1},\, x_{k+1} - x_k\rangle + \langle x_k - x_{k-1},\, x_k - x^\star\rangle\\ &= \langle x_k - x_{k-1},\, x_{k+1} - x_k\rangle + \big(E_{x,k} + \varphi_k - \varphi_{k-1}\big), \end{aligned} \tag{A.7}$$

where we applied the usual Pythagoras relation to $\langle x_k - x_{k-1},\, x_k - x^\star\rangle$,

$$2\langle c_1 - c_2,\, c_1 - c_3\rangle = \|c_1 - c_2\|^2 + \|c_1 - c_3\|^2 - \|c_2 - c_3\|^2.$$

Putting (A.7) back into (A.6) yields

$$\begin{aligned} &\varphi_{k+1} - \varphi_k - a_k(\varphi_k - \varphi_{k-1})\\ &\le -E_{x,k+1} - \gamma_k\langle B(y_{b,k}) - B(x^\star),\, x_{k+1} - x^\star\rangle + a_k\langle x_k - x_{k-1},\, x_{k+1} - x_k\rangle + a_k E_{x,k}. \end{aligned} \tag{A.8}$$

Since $B$ is $\beta$-cocoercive, then

$$\begin{aligned}\langle B(y_{b,k}) - B(x^\star),\, x_{k+1} - x^\star\rangle &= \langle B(y_{b,k}) - B(x^\star),\, x_{k+1} - y_{b,k} + y_{b,k} - x^\star\rangle \\ &\geq \beta\|B(y_{b,k}) - B(x^\star)\|^2 + \langle B(y_{b,k}) - B(x^\star),\, x_{k+1} - y_{b,k}\rangle \\ &\geq \beta\|B(y_{b,k}) - B(x^\star)\|^2 - \beta\|B(y_{b,k}) - B(x^\star)\|^2 - \frac{1}{2\beta}E_{b,k+1} \\ &= -\frac{1}{2\beta}E_{b,k+1}.\end{aligned} \tag{A.9}$$

Denote $\mu_k = 1 - \frac{\gamma_k}{2\beta} \in [\frac{\bar\epsilon}{2\beta}, 1 - \frac{\epsilon}{2\beta}]$, $\nu_k = a_k - \frac{\gamma_k b_k}{2\beta}$ and $v_k = x_{k+1} - x_k - \frac{\nu_k}{\mu_k}(x_k - x_{k-1})$. Substituting (A.9) back into (A.8), and since $E_{b,k+1} = E_{x,k+1} + b_k^2 E_{x,k} + b_k\langle x_k - x_{k+1},\, x_k - x_{k-1}\rangle$, we get

$$\begin{aligned}&\varphi_{k+1} - \varphi_k - a_k(\varphi_k - \varphi_{k-1}) \\ &\leq -E_{x,k+1} + \frac{\gamma_k}{2\beta}E_{b,k+1} + a_k\langle x_k - x_{k-1},\, x_{k+1} - x_k\rangle + a_k E_{x,k} \\ &= -\mu_k E_{x,k+1} + \big(a_k - \frac{\gamma_k b_k}{2\beta}\big)\langle x_k - x_{k-1},\, x_{k+1} - x_k\rangle + \big(a_k + \frac{\gamma_k b_k^2}{2\beta}\big)E_{x,k} \\ &= -\frac{\mu_k}{2}\|x_k - x_{k+1}\|^2 + \nu_k\langle x_k - x_{k-1},\, x_{k+1} - x_k\rangle + \big(a_k + \frac{\gamma_k b_k^2}{2\beta}\big)E_{x,k} \\ &= \big(-\frac{\mu_k}{2}\|x_{k+1} - x_k - \frac{\nu_k}{\mu_k}(x_k - x_{k-1})\|^2 + \frac{\nu_k^2}{\mu_k}E_{x,k}\big) + \big(a_k + \frac{\gamma_k b_k^2}{2\beta}\big)E_{x,k} \\ &= -\frac{\mu_k}{2}\|v_k\|^2 + \big(a_k + \frac{\nu_k^2}{\mu_k} + \frac{\gamma_k b_k^2}{2\beta}\big)E_{x,k} \leq -\frac{\mu_k}{2}\|v_k\|^2 + \big(\frac{2a_k}{\mu_k} + \frac{\gamma_k b_k}{2\beta}\big)E_{x,k} \\ &\leq -\frac{\mu_k}{2}\|v_k\|^2 + \big(\frac{4\beta}{\bar\epsilon}a_k + (1 - \frac{\bar\epsilon}{2\beta})b_k\big)E_{x,k}.\end{aligned} \tag{A.10}$$

Denote $\theta_k = \varphi_k - \varphi_{k-1}$ and $\delta_k = \big(\frac{4\beta}{\bar\epsilon}a_k + (1 - \frac{\bar\epsilon}{2\beta})b_k\big)E_{x,k}$. We then arrive at the following key estimate

$$\theta_{k+1} \leq -\frac{\mu_k}{2}\|v_k\|^2 + a_k\theta_k + \delta_k. \tag{A.11}$$

If $a_k \in ]0, \bar a]$, (A.11) yields

$$\theta_{k+1} \leq -\frac{\mu_k}{2}\|v_k\|^2 + a_k\theta_k + \delta_k \leq a_k\theta_k + \delta_k \leq a_k[\theta_k]_+ + \delta_k, \tag{A.12}$$

where $[\theta]_+ = \max\{\theta, 0\}$. As a result, we have

$$[\theta_{k+1}]_+ \leq \bar a[\theta_k]_+ + \delta_k.$$

Assumption (2.1) is equivalent to the fact that $\delta_k$ is summable. Therefore, using that $\bar a < 1$ and applying [21, Lemma 3.1(iv)], it follows that $[\theta_k]_+$ is summable. Therefore,

$$\varphi_{k+1} - \sum_{j=1}^{k+1}[\theta_j]_+ \leq \varphi_{k+1} - \theta_{k+1} - \sum_{j=1}^{k}[\theta_j]_+ = \varphi_k - \sum_{j=1}^{k}[\theta_j]_+.$$

It follows that the sequence $(\varphi_k - \sum_{j=1}^{k}[\theta_j]_+)_{k\in\mathbb{N}}$ is decreasing and bounded below, hence convergent, whence we deduce that $\varphi_k$ is also convergent.

If $a_k \equiv 0$, (A.10) entails

$$\varphi_{k+1} \leq \varphi_k + \delta_k.$$

We then conclude that the sequence $(x_k)_{k\in\mathbb{N}}$ is quasi-Fejér monotone (of type III) relative to $\mathrm{zer}(A+B)$ [21, Definition 1.1(3)], and thus $\varphi_k$ is convergent [21, Proposition 3.6].

In summary, for $a_k \in [0,\bar{a}]$, $\lim_{k\to\infty}\|x_k - x^\star\|$ exists for any $x^\star \in \mathrm{zer}(A+B)$.

By assumption (2.1), $a_k(x_k - x_{k-1}) \to 0$ and $b_k(x_k - x_{k-1}) \to 0$, and thus

$$\tfrac{\nu_k}{\mu_k}(x_k - x_{k-1}) \to 0, \tag{A.13}$$

since $\mu_k \geq \frac{\bar{\epsilon}}{2\beta} > 0$. Moreover, from (A.12), we obtain

$$\textstyle\sum_{k\in\mathbb{N}}\|v_k\|^2 \leq \frac{4\beta}{\bar{\epsilon}}\big(\bar{a}\varphi_0 + \sum_{k\in\mathbb{N}}(\bar{a}[\theta_k]_+ + \delta_k)\big) < +\infty.$$

Consequently, $v_k \to 0$. Combining this with (A.13), we get that $x_{k+1} - x_k \to 0$. In turn, $y_{a,k} - x_{k+1} \to 0$ and $y_{b,k} - x_{k+1} \to 0$.

Let $\bar{x}$ be a weak cluster point of $(x_k)_{k\in\mathbb{N}}$, and let us fix a subsequence, say $x_{k_j} \rightharpoonup \bar{x}$. We get from (1.4) that

$$u_{k_j} \stackrel{\text{def}}{=} \frac{y_{a,k_j} - x_{k_j+1}}{\gamma_{k_j}} - B(y_{b,k_j}) \in A(x_{k_j+1}).$$

Since $B$ is cocoercive and $y_{b,k_j} \rightharpoonup \bar{x}$, we have $B(y_{b,k_j}) \to B(\bar{x})$. In turn, $u_{k_j} \to -B(\bar{x})$ since $\gamma_k \geq \underline{\epsilon} > 0$. Since $(x_{k_j+1}, u_{k_j}) \in \mathrm{gph}\,A$, and the graph of the maximal monotone operator $A$ is sequentially weakly-strongly closed in $\mathcal{H}\times\mathcal{H}$, we get that $-B(\bar{x}) \in A(\bar{x})$, *i.e.* $\bar{x}$ is a solution of $(\mathcal{P}_{\text{inc}})$. Opial's Theorem [49] then concludes the proof. □

**Proof of Theorem 2.3.** From (A.10), we apply Young's inequality to get

$$\begin{aligned}
&\varphi_{k+1} - \varphi_k - a_k(\varphi_k - \varphi_{k-1})\\
&\leq \Big(\frac{\gamma_k}{2\beta} - 1\Big)E_{x,k+1} + \Big(a_k - \frac{\gamma_k b_k}{2\beta}\Big)\langle x_k - x_{k-1},\, x_{k+1} - x_k\rangle + \Big(a_k + \frac{\gamma_k b_k^2}{2\beta}\Big)E_{x,k}\\
&\leq \Big(\frac{\gamma_k}{2\beta} - 1\Big)E_{x,k+1} + \Big|a_k - \frac{\gamma_k b_k}{2\beta}\Big|\frac{1}{2}\big(\|x_{k+1} - x_k\|^2 + \|x_k - x_{k-1}\|^2\big) + \Big(\frac{\gamma_k}{2\beta}b_k^2 + a_k\Big)E_{x,k}\\
&= \Big(\frac{\gamma_k}{2\beta} - 1 + \Big|a_k - \frac{\gamma_k b_k}{2\beta}\Big|\Big)E_{x,k+1} + \Big(\frac{\gamma_k}{2\beta}b_k^2 + a_k + \Big|a_k - \frac{\gamma_k b_k}{2\beta}\Big|\Big)E_{x,k}\\
&= S_k E_{x,k+1} + T_k E_{x,k},
\end{aligned}$$

where $S_k = \frac{\gamma_k}{2\beta} - 1 + |a_k - \frac{\gamma_k b_k}{2\beta}|$, $T_k = \frac{\gamma_k}{2\beta}b_k^2 + a_k + |a_k - \frac{\gamma_k b_k}{2\beta}|$. Suppose $a_k$, $b_k$ and $\gamma_k$ are non-decreasing such that $S_k$, $T_k$ are also non-decreasing. Define $\phi_k = \varphi_k - a_k\varphi_{k-1} + T_k E_{x,k}$, then

$$\begin{aligned}
\phi_{k+1} - \phi_k &= (\varphi_{k+1} - a_{k+1}\varphi_k + T_{k+1}E_{x,k+1}) - (\varphi_k - a_k\varphi_{k-1} + T_k E_{x,k})\\
&\leq (\varphi_{k+1} - \varphi_k) - a_k(\varphi_k - \varphi_{k-1}) + T_{k+1}E_{x,k+1} - T_k E_{x,k}\\
&\leq S_k E_{x,k+1} + T_k E_{x,k} + T_{k+1}E_{x,k+1} - T_k E_{x,k} = (S_k + T_{k+1})E_{x,k+1}.
\end{aligned} \tag{A.14}$$

Case 1) $a_k \in [0,\bar{a}]$, $b_k \in [0,\bar{b}]$, $b_k \leq a_k$. We have $\frac{\gamma_k}{2\beta}b_k < a_k$, then from (A.14), and under the second condition in (2.4),

$$\phi_{k+1} - \phi_k \leq (S_{k+1} + T_{k+1})E_{x,k+1} = \Big((3a_{k+1} - 1) + \frac{\gamma_{k+1}}{2\beta}(1 - b_{k+1})^2\Big)E_{x,k+1} \leq -\tau E_{x,k+1}, \tag{A.15}$$

Case 2) $a_k \in [0, \bar{a}]$, $b_k \in ]0, \bar{b}]$, $a_k < b_k$. Since $S_k, T_k$ are non-decreasing, then from (A.14) we have,

$$\begin{aligned}
&\phi_{k+1} - \phi_k \le (S_{k+1} + T_{k+1})E_{x,k+1} \\
&\le \Big(\tfrac{\gamma_{k+1}}{2\beta} - 1 + |a_{k+1} - \tfrac{\gamma_{k+1}}{2\beta} b_{k+1}| + \tfrac{\gamma_{k+1}}{2\beta} b_{k+1}^2 + a_{k+1} + |a_{k+1} - \tfrac{\gamma_{k+1}}{2\beta} b_{k+1}|\Big) E_{x,k+1}.
\end{aligned}$$

Next we discuss the relationship between $a_{k+1}$ and $\frac{\gamma_{k+1}}{2\beta} b_{k+1}$, which splits into two subcases.

(i) If $\frac{\gamma_{k+1}}{2\beta} b_{k+1} \le a_{k+1}$, $k \in \mathbb{N}$, then from the second condition in (2.4),

$$\begin{aligned}
\phi_{k+1} - \phi_k &\le \Big(\tfrac{\gamma_{k+1}}{2\beta} - 1 + a_{k+1} - \tfrac{\gamma_{k+1}}{2\beta} b_{k+1} + \tfrac{\gamma_{k+1}}{2\beta} b_{k+1}^2 + 2a_{k+1} - \tfrac{\gamma_{k+1}}{2\beta} b_{k+1}\Big) E_{x,k+1} \\
&= \Big((3a_{k+1} - 1) + \tfrac{\gamma_{k+1}}{2\beta}(1 - b_{k+1})^2\Big) E_{x,k+1} \le -\tau E_{x,k+1}.
\end{aligned} \tag{A.16}$$

(ii) If $a_{k+1} < \frac{\gamma_{k+1}}{2\beta} b_{k+1}$, $k \in \mathbb{N}$, then from the first condition of (2.4), we have

$$\begin{aligned}
\phi_{k+1} - \phi_k &\le \Big(\tfrac{\gamma_{k+1}}{2\beta} - 1 + \tfrac{\gamma_{k+1}}{2\beta} b_{k+1} - a_{k+1} + \tfrac{\gamma_{k+1}}{2\beta} b_{k+1}^2 + \tfrac{\gamma_{k+1}}{2\beta} b_{k+1}\Big) E_{x,k+1} \\
&= \Big(-(1 + a_{k+1}) + \tfrac{\gamma_{k+1}}{2\beta}(1 + b_{k+1})^2\Big) E_{x,k+1} \le -\tau E_{x,k+1}.
\end{aligned} \tag{A.17}$$

From (A.15) (respectively (A.16) or (A.17)), $\phi_k$ is non-increasing. Therefore, we have

$$\varphi_k - \bar{a}\varphi_{k-1} \le \phi_k \le \phi_1 \Longrightarrow \varphi_k \le \bar{a}^k \varphi_0 + \phi_1 \textstyle\sum_{j=0}^{k-1} \bar{a}^j \le \bar{a}^k \varphi_0 + \dfrac{\phi_1}{1 - \bar{a}}.$$

In the meanwhile, from (A.15) we have

$$\begin{aligned}
&\phi_{k+1} - \phi_1 \le -\tau \textstyle\sum_{j=0}^{k} E_{x,j+1} \\
\Longrightarrow\ &\textstyle\sum_{j=0}^{k} E_{x,j} \le \dfrac{1}{\tau}(\phi_1 - \phi_{k+1}) \le \dfrac{1}{\tau}(\phi_1 + \bar{a}\varphi_k) \le \dfrac{1}{\tau}\Big(\bar{a}^{k+1}\varphi_0 + \dfrac{\phi_1}{1 - \bar{a}}\Big) < +\infty,
\end{aligned}$$

which means that the summability condition in (2.2) is satisfied. The rest of the proof follows the same arguments as in those in the last part of the proof of Theorem 2.1. □

Denote $A^{\varepsilon}$ the $\varepsilon$-enlargements of $A$.

**Proof of Proposition 2.5.** Let $x^\star \in \operatorname{zer}(A + B)$. Recall from (2.5) that

$$y_{a,k} - \gamma_k B(y_{b,k}) - \gamma_k \xi_k - x_{k+1} \in \gamma_k A^{\varepsilon_k}(x_{k+1}).$$

Thus, we get

$$\langle y_{a,k} - x_{k+1} - \gamma_k (B(y_{b,k}) - B(x^\star)) - \gamma_k \xi_k,\ x_{k+1} - x^\star \rangle \ge -\gamma_k \varepsilon_k.$$

Combining this with (A.5), we obtain

$$\varphi_k - \varphi_{k+1} \ge E_{x,k+1} + \gamma_k \langle B(y_{b,k}) - B(x^\star) + \xi_k,\ x_{k+1} - x^\star \rangle - a_k \langle x_k - x_{k-1},\ x_{k+1} - x^\star \rangle - \gamma_k \varepsilon_k.$$

Continuing as after (A.6) in the proof of Theorem 2.1, we obtain the key estimate

$$\begin{aligned}
\theta_{k+1} &\le -\tfrac{\mu_k}{2}\|v_k\|^2 + a_k \theta_k + \delta_k + \gamma_k \varepsilon_k + \gamma_k \langle \xi_k,\ x_{k+1} - x^\star \rangle \\
&\le -\tfrac{\mu_k}{2}\|v_k\|^2 + \bar{a}\theta_k + \delta_k + \overline{\gamma}\varepsilon_k + \sqrt{2}\overline{\gamma}\|\xi_k\|\sqrt{\varphi_{k+1}},
\end{aligned} \tag{A.18}$$

where $\overline{\gamma} = (2\beta - \bar{\epsilon})$, $\theta_k$, $\delta_k$ and $v_k$ are as defined in (A.11). This yields

$$\theta_{k+1} \le \bar{a}^k \theta_1 + \textstyle\sum_{j=1}^{k} \bar{a}^{k-j} \big(\delta_j + \overline{\gamma}\varepsilon_j + \sqrt{2}\overline{\gamma}\|\xi_j\|\sqrt{\varphi_{j+1}}\big).$$

(i) $a_k \in ]0, \bar{a}]$: summing up the last inequality, we get

$$\sum_{m=1}^{k} \theta_{m+1} = \varphi_{k+1} - \varphi_1 \leq \frac{1}{1-\bar{a}}\big(\varphi_1 - \varphi_0 + \sum_{m\in\mathbb{N}} \delta_m + \bar{\gamma} \sum_{m\in\mathbb{N}} \varepsilon_{m+1}\big) + \sqrt{2\bar{\gamma}} \sum_{m=1}^{k} m\|\xi_m\|\sqrt{\varphi_{m+1}},$$

which entails

$$\varphi_{k+1} \leq c + \sqrt{2\bar{\gamma}} \sum_{m=1}^{k} m\|\xi_m\|\sqrt{\varphi_m} \leq c + \sqrt{2\bar{\gamma}} \sum_{m=1}^{k+1} m\|\xi_m\|\sqrt{\varphi_{m+1}}, \tag{A.19}$$

where $c = \varphi_1 + \frac{1}{1-\bar{a}}(\varphi_1 + \sum_{m\in\mathbb{N}} \delta_m + \bar{\gamma} \sum_{m\in\mathbb{N}} \varepsilon_m) \geq 0$. By assumption on the sequences $(\varepsilon_m)_{m\in\mathbb{N}}$ and $(\delta_m)_{m\in\mathbb{N}}$, $c$ is bounded. Using the fact that $(m\|\xi_m\|)_{m\in\mathbb{N}}$ is summable, it can be easily shown, *e.g.* [6, Lemma A.9], that since the sequence $(\varphi_k)_{k\in\mathbb{N}}$ satisfies (A.19), it also obeys $\varphi_k \leq \sqrt{c} + \sum_{j\in\mathbb{N}} j\|\xi_j\| < +\infty$. Denote $t = \sqrt{c} + \sum_{j\in\mathbb{N}} j\|\xi_j\|$. (A.18) then becomes

$$\theta_{k+1} \leq -\frac{\mu_k}{2}\|v_k\|^2 + \bar{a}\theta_k + \delta_k + \bar{\gamma}\varepsilon_k + \sqrt{2t\bar{\gamma}}\|\xi_k\|,$$

which is of the form (A.12), where $\delta_k$ is replaced by $\delta_k + \bar{\gamma}\varepsilon_k + \sqrt{2\bar{\gamma}}\sqrt{t}\|\xi_k\|$, and the latter is a summable sequence. With the same arguments as those after (A.12) for $a_k \in ]0, \bar{a}]$, we deduce that $\varphi_k$ is convergent.

(ii) $a_k \equiv 0$: in this case, (A.18) reduces to

$$\varphi_{k+1} \leq \varphi_k + \delta_k + \bar{\gamma}\varepsilon_k + \sqrt{2\bar{\gamma}}\|\xi_k\|\sqrt{\varphi_{k+1}} \leq \varphi_1 + \sum_{j\in\mathbb{N}} \delta_j + \bar{\gamma} \sum_{j\in\mathbb{N}} \varepsilon_j + \sqrt{2\bar{\gamma}} \sum_{j=1}^{k+1} \|\xi_j\|\sqrt{\varphi_{j+1}}.$$

Again, by virtue of [6, Lemma A.9] and summability of the sequences $(\delta_j)_{j\in\mathbb{N}}$, $(\varepsilon_j)_{j\in\mathbb{N}}$ and $(\|\xi_j\|)_{j\in\mathbb{N}}$, we have $\varphi_k \leq t = \sqrt{\varphi_1 + \sum_{j\in\mathbb{N}}(\delta_j + \bar{\gamma}\varepsilon_j + \|\xi_j\|)} < +\infty$. Consequently, we have

$$\varphi_{k+1} \leq \varphi_k + \delta_k + \bar{\gamma}\varepsilon_k + \sqrt{2t\bar{\gamma}}\|\xi_k\|.$$

That is, the sequence $(x_k)_{k\in\mathbb{N}}$ is quasi-Fejér monotone (of type III) relative to $\mathrm{zer}(A+B)$, and thus $\varphi_k$ is convergent.

In summary, for $a_k \in [0, \bar{a}]$, $\lim_{k\to\infty} \|x_k - x^\star\|$ exists for any $x^\star \in \mathrm{zer}(A+B)$.

The rest of the proof is patterned after the last part of the proof of Theorem 2.1, where we now use the fact that $\xi_k \to 0$ by assumption, and that the graph of $A^{\cdot} : \mathbb{R}_+ \times \mathcal{H} \rightrightarrows \mathcal{H}$ is weakly-strongly sequentially closed in $\mathbb{R}_+ \times \mathcal{H} \times \mathcal{H}$ [57, Proposition 3.4(b)]. □

# B Proofs of Section 4

## B.1 Riemannian Geometry

Let $\mathcal{M}$ be a $C^2$-smooth embedded submanifold of $\mathbb{R}^n$ around a point $x$. With some abuse of terminology, we shall state $C^2$-manifold instead of $C^2$-smooth embedded submanifold of $\mathbb{R}^n$. The natural embedding of a submanifold $\mathcal{M}$ into $\mathbb{R}^n$ permits to define a Riemannian structure and to introduce geodesics on $\mathcal{M}$, and we simply say $\mathcal{M}$ is a Riemannian manifold. We denote respectively $\mathcal{T}_{\mathcal{M}}(x)$ and $\mathcal{N}_{\mathcal{M}}(x)$ the tangent and normal space of $\mathcal{M}$ at point near $x$ in $\mathcal{M}$.

**Exponential map** Geodesics generalize the concept of straight lines in $\mathbb{R}^n$, preserving the zero acceleration characteristic, to manifolds. Roughly speaking, a geodesic is locally the shortest path between two points on $\mathcal{M}$. We denote by $\mathfrak{g}(t; x, h)$ the value at $t \in \mathbb{R}$ of the geodesic starting at $\mathfrak{g}(0; x, h) = x \in \mathcal{M}$ with velocity $\dot{\mathfrak{g}}(t; x, h) = \frac{d\mathfrak{g}}{dt}(t; x, h) = h \in \mathcal{T}_{\mathcal{M}}(x)$ (which is uniquely defined). For every $h \in \mathcal{T}_{\mathcal{M}}(x)$, there exists an interval $I$ around 0 and a unique geodesic $\mathfrak{g}(t; x, h) : I \to \mathcal{M}$ such that $\mathfrak{g}(0; x, h) = x$ and $\dot{\mathfrak{g}}(0; x, h) = h$. The mapping

$$\mathrm{Exp}_x : \mathcal{T}_{\mathcal{M}}(x) \to \mathcal{M}, \ h \mapsto \mathrm{Exp}_x(h) = \mathfrak{g}(1; x, h),$$

is called *Exponential map*. Given $x, x' \in \mathcal{M}$, the direction $h \in \mathcal{T}_{\mathcal{M}}(x)$ we are interested in is such that

$$\mathrm{Exp}_x(h) = x' = \mathfrak{g}(1; x, h).$$

**Parallel translation** Given two points $x, x' \in \mathcal{M}$, let $\mathcal{T}_{\mathcal{M}}(x), \mathcal{T}_{\mathcal{M}}(x')$ be their corresponding tangent spaces. Define

$$\tau : \mathcal{T}_{\mathcal{M}}(x) \to \mathcal{T}_{\mathcal{M}}(x'),$$

the parallel translation along the unique geodesic joining $x$ to $x'$, which is isomorphism and isometry w.r.t. the Riemannian metric.

**Riemannian gradient and Hessian** For a vector $v \in \mathcal{N}_{\mathcal{M}}(x)$, the Weingarten map of $\mathcal{M}$ at $x$ is the operator $\mathfrak{W}_x(\cdot, v) : \mathcal{T}_{\mathcal{M}}(x) \to \mathcal{T}_{\mathcal{M}}(x)$ defined by

$$\mathfrak{W}_x(\cdot, v) = -\mathrm{P}_{\mathcal{T}_{\mathcal{M}}(x)} \mathrm{d}V[h],$$

where $V$ is any local extension of $v$ to a normal vector field on $\mathcal{M}$. The definition is independent of the choice of the extension $V$, and $\mathfrak{W}_x(\cdot, v)$ is a symmetric linear operator which is closely tied to the second fundamental form of $\mathcal{M}$, see [20, Proposition II.2.1].

Let $G$ be a real-valued function which is $C^2$ along the $\mathcal{M}$ around $x$. The covariant gradient of $G$ at $x' \in \mathcal{M}$ is the vector $\nabla_{\mathcal{M}} G(x') \in \mathcal{T}_{\mathcal{M}}(x')$ defined by

$$\langle \nabla_{\mathcal{M}} G(x'),\, h \rangle = \frac{d}{dt} G\big(\mathrm{P}_{\mathcal{M}}(x' + th)\big)\big|_{t=0}, \ \forall h \in \mathcal{T}_{\mathcal{M}}(x'),$$

where $\mathrm{P}_{\mathcal{M}}$ is the projection operator onto $\mathcal{M}$. The covariant Hessian of $G$ at $x'$ is the symmetric linear mapping $\nabla^2_{\mathcal{M}} G(x')$ from $\mathcal{T}_{\mathcal{M}}(x')$ to itself which is defined as

$$\langle \nabla^2_{\mathcal{M}} G(x')h,\, h \rangle = \frac{d^2}{dt^2} G\big(\mathrm{P}_{\mathcal{M}}(x' + th)\big)\big|_{t=0}, \ \forall h \in \mathcal{T}_{\mathcal{M}}(x'). \tag{B.1}$$

This definition agrees with the usual definition using geodesics or connections [42]. Now assume that $\mathcal{M}$ is a Riemannian embedded submanifold of $\mathbb{R}^n$, and that a function $G$ has a $C^2$-smooth restriction on $\mathcal{M}$. This can be characterized by the existence of a $C^2$-smooth extension (representative) of $G$, *i.e.* a $C^2$-smooth function $\widetilde{G}$ on $\mathbb{R}^n$ such that $\widetilde{G}$ agrees with $G$ on $\mathcal{M}$. Thus, the Riemannian gradient $\nabla_{\mathcal{M}} G(x')$ is also given by

$$\nabla_{\mathcal{M}} G(x') = \mathrm{P}_{\mathcal{T}_{\mathcal{M}}(x')} \nabla \widetilde{G}(x'), \tag{B.2}$$

and $\forall h \in \mathcal{T}_{\mathcal{M}}(x')$, the Riemannian Hessian reads

$$\begin{aligned} \nabla^2_{\mathcal{M}} G(x')h &= \mathrm{P}_{\mathcal{T}_{\mathcal{M}}(x')} \mathrm{d}(\nabla_{\mathcal{M}} G)(x')[h] = \mathrm{P}_{\mathcal{T}_{\mathcal{M}}(x')} \mathrm{d}\big(x' \mapsto \mathrm{P}_{\mathcal{T}_{\mathcal{M}}(x')} \nabla_{\mathcal{M}} \widetilde{G}\big)[h] \\ &= \mathrm{P}_{\mathcal{T}_{\mathcal{M}}(x')} \nabla^2 \widetilde{G}(x')h + \mathfrak{W}_{x'}\big(h, \mathrm{P}_{\mathcal{N}_{\mathcal{M}}(x')} \nabla \widetilde{G}(x')\big), \end{aligned} \tag{B.3}$$

where the last equality comes from [2, Theorem 1]. When $\mathcal{M}$ is an affine or linear subspace of $\mathbb{R}^n$, then obviously $\mathcal{M} = x + \mathcal{T}_{\mathcal{M}}(x)$, and $\mathfrak{W}_{x'}(h, \mathrm{P}_{\mathcal{N}_{\mathcal{M}}(x')}\nabla \widetilde{G}(x')) = 0$, hence (B.3) reduces to

$$\nabla^2_{\mathcal{M}} G(x') = \mathrm{P}_{\mathcal{T}_{\mathcal{M}}(x')} \nabla^2 \widetilde{G}(x') \mathrm{P}_{\mathcal{T}_{\mathcal{M}}(x')}.$$

See [35, 20] for more materials on differential and Riemannian manifolds.

The following lemmas summarize two key properties that we will need throughout.

**Lemma B.1.** *Let $x \in \mathcal{M}$, and $x_k$ a sequence converging to $x$ in $\mathcal{M}$. Denote $\tau_k : \mathcal{T}_{\mathcal{M}}(x) \to \mathcal{T}_{\mathcal{M}}(x_k)$ be the parallel translation along the unique geodesic joining $x$ to $x_k$. Then, for any bounded vector $u \in \mathbb{R}^n$, we have*

$$(\tau_k^{-1} \mathrm{P}_{\mathcal{T}_{\mathcal{M}}(x_k)} - \mathrm{P}_{\mathcal{T}_{\mathcal{M}}(x)}) u = o(\|u\|).$$

**Proof.** From [1, Chapter 5], we deduce that for $k$ sufficiently large,

$$\tau_k^{-1} = \mathrm{P}_{\mathcal{T}_{\mathcal{M}}(x)} + o(\|x_k - x\|).$$

In addition, locally near $x$ along $\mathcal{M}$, the operator $x \mapsto \mathrm{P}_{\mathcal{T}_{\mathcal{M}}(x)}$ is $C^1$, hence,

$$\begin{aligned} \lim_{k\to\infty} \frac{\|(\tau_k^{-1} \mathrm{P}_{\mathcal{T}_{\mathcal{M}}(x_k)} - \mathrm{P}_{\mathcal{T}_{\mathcal{M}}(x)})u\|}{\|u\|} &\le \lim_{k\to\infty} \frac{\|\mathrm{P}_{\mathcal{T}_{\mathcal{M}}(x)}(\mathrm{P}_{\mathcal{T}_{\mathcal{M}}(x_k)} - \mathrm{P}_{\mathcal{T}_{\mathcal{M}}(x)})\|\|u\|}{\|u\|} + o(\|x_k - x\|) \\ &\le \lim_{k\to\infty} \|\mathrm{P}_{\mathcal{T}_{\mathcal{M}}(x_k)} - \mathrm{P}_{\mathcal{T}_{\mathcal{M}}(x)}\| + o(\|x_k - x\|) = 0. \end{aligned}$$

□

**Lemma B.2.** *Let $x, x'$ be two close points in $\mathcal{M}$, denote $\tau : \mathcal{T}_{\mathcal{M}}(x) \to \mathcal{T}_{\mathcal{M}}(x')$ the parallel translation along the unique geodesic joining $x$ to $x'$. The Riemannian Taylor expansion of $\Phi \in C^2(\mathcal{M})$ around $x$ reads,*

$$\tau^{-1} \nabla_{\mathcal{M}} \Phi(x') = \nabla_{\mathcal{M}} \Phi(x) + \nabla^2_{\mathcal{M}} \Phi(x) \mathrm{P}_{\mathcal{T}_{\mathcal{M}}(x)}(x' - x) + o(\|x' - x\|).$$

**Proof.** Since $x, x' \in \mathcal{M}$ are close, we have $x' = \mathrm{Exp}_x(h)$ for some $h \in \mathcal{T}_{\mathcal{M}}(x)$ small enough, and thus, the Taylor expansion [56, Remark 4.2] of $\nabla_{\mathcal{M}} \Phi$ around $x$ reads

$$\tau^{-1} \nabla_{\mathcal{M}} \Phi(x') = \nabla_{\mathcal{M}} \Phi(x) + \nabla^2_{\mathcal{M}} \Phi(x) h + o(\|h\|). \tag{B.4}$$

Moreover, form the proof of [42, Theorem 4.9], one can show that

$$\mathrm{P}_{\mathcal{T}_{\mathcal{M}}(x)}(x') = \mathrm{P}_{\mathcal{T}_{\mathcal{M}}(x)}(\mathrm{Exp}_x(h)) = \mathrm{P}_{\mathcal{T}_{\mathcal{M}}(x)}(x) + h + o(\|h\|^2).$$

Substituting back into (B.4) we get the claimed result. □

## B.2 Proofs

**Proof of Proposition 4.1.**

(i) Since $F$ is locally $C^2$ around $x^\star$, there exists $\epsilon > 0$ sufficiently small such that for any $\delta \in \mathbb{B}_\epsilon(0)$, we have

$$\begin{aligned} \Phi(x^\star + \delta) - \Phi(x^\star) &= F(x^\star + \delta) - F(x^\star) - \langle \nabla F(x^\star), \delta \rangle + R(x^\star + \delta) - R(x^\star) + \langle \nabla F(x^\star), \delta \rangle \\ &= \frac{1}{2} \langle \delta, \nabla^2 F(x^\star + t\delta)\delta \rangle + R(x^\star + \delta) - R(x^\star) + \langle \nabla F(x^\star), \delta \rangle, \quad t \in ]0, 1[. \end{aligned}$$

Let $x_t = x^\star + t\delta \in \mathbb{B}_\epsilon(x^\star)$. Since (RI) holds and $\nabla^2 F(x)$ depends continuously on $x \in \mathbb{B}_\epsilon(x^\star)$, we have $\mathrm{P}_{T_{x^\star}} \nabla^2 F(x) \mathrm{P}_{T_{x^\star}} \succeq \alpha \mathrm{Id}$ for any such $x$. This holds in particular at $x_t$. We then distinguish two cases.

(a) $\delta \notin \ker(\nabla^2 F(x_t))$. In this case, it is clear that

$$\Phi(x^\star + \delta) - \Phi(x^\star) \geq \frac{1}{2}\langle \delta, \nabla^2 F(x_t)\delta \rangle \geq \alpha/2\|\delta\|^2 > 0$$

since $F$ is convex and locally $C^2$, and $R$ is convex with $-\nabla F(x^\star) \in \partial R(x^\star)$.

(b) $\delta \in \ker(\nabla^2 F(x_t)) \setminus \{0\}$. Since $R$ is a proper closed convex function, it is sub-differentially regular at $x^\star$. Moreover $\partial R(x^\star) \neq \emptyset$ ($-\nabla F(x^\star)$ is in it), and thus the directional derivative $R'(x^\star, \cdot)$ is proper and closed, and it is the support of $\partial R(x^\star)$ [54, Theorem 8.30]. It then follows from the separation theorem [32, Theorem V.2.2.3] that

$$-\nabla F(x^\star) \in \mathrm{ri}(\partial R(x^\star)) \iff R'(x^\star, \delta) > -\langle \nabla F(x^\star), \delta \rangle, \forall \delta \text{ s.t. } R'(x^\star; \delta) + R'(x^\star; -\delta) > 0.$$

As $\ker(R'(x^\star; \cdot)) = T_{x^\star}$ [61, Proposition 3(iii) and Lemma 10], and in view of (RI), we get

$$\begin{aligned} -\nabla F(x^\star) \in \mathrm{ri}(\partial R(x^\star)) &\iff R'(x^\star; \delta) > -\langle \nabla F(x^\star), \delta \rangle, \forall \delta \notin T_{x^\star} \\ &\implies R'(x^\star; \delta) > -\langle \nabla F(x^\star), \delta \rangle, \forall \delta \in \ker(\nabla^2 F(x_t)) \setminus \{0\}. \end{aligned}$$

Combining this with classical properties of the directional derivative of a convex function yields

$$\begin{aligned} \Phi(x^\star + \delta) - \Phi(x^\star) &= R(x^\star + \delta) - R(x^\star) + \langle \nabla F(x^\star), \delta \rangle \\ &\geq R'(x^\star; \delta) + \langle \nabla F(x^\star), \delta \rangle > 0, \end{aligned}$$

which concludes the first claim.

(ii) Let $\Psi$ as defined in the proof of Lemma 4.3. If $R \in \mathrm{PSF}_{x^\star}(\mathcal{M}_{x^\star})$, the Riemannian Hessian of $\Phi$ reads

$$\nabla^2_{\mathcal{M}_{x^\star}} \Phi(x^\star) = \mathrm{P}_{T_{x^\star}} \nabla F(x^\star) \mathrm{P}_{T_{x^\star}} + \nabla^2_{\mathcal{M}_{x^\star}} \Psi(x^\star).$$

In view of Lemma 4.3(i), $\nabla^2_{\mathcal{M}_{x^\star}} \Psi(x^\star)$ is positive semi-definite on $T_{x^\star}$. On the other hand, hypothesis (RI) entails positive definiteness of $\mathrm{P}_{T_{x^\star}} \nabla F(x^\star) \mathrm{P}_{T_{x^\star}}$. Altogether, this shows that $\nabla^2_{\mathcal{M}_{x^\star}} \Phi(x^\star)$ is positive definite on $T_{x^\star} \setminus \{0\}$. Local quadratic growth of $\Phi$ near $x^\star$ then follows by combining [37, Definition 5.4], [42, Theorem 3.4] and [30, Theorem 6.2].

□

**Proof of Lemma 4.3.** By definition of $U$, $Uh = 0$ for any $h \in T^\perp_{x^\star}$. Thus, in the following we only examine the case $h \in T_{x^\star}$.

(i) Let $\Psi(x) \stackrel{\text{def}}{=} R(x) + \langle x, \nabla F(x^\star) \rangle$. From the smooth perturbation rule of partial smoothness [37, Corollary 4.7], $\Psi \in \mathrm{PSF}_{x^\star}(\mathcal{M}_{x^\star})$. Moreover, from Fact 3.3 and normal sharpness, the Riemannian Hessian of $\Psi$ at $x^\star$ is such that, $\forall h \in T_{x^\star}$,

$$\begin{aligned} \gamma \nabla^2_{\mathcal{M}_{x^\star}} \Psi(x^\star) h &= \gamma \mathrm{P}_{T_{x^\star}} \nabla^2 \widetilde{\Psi}(x^\star) h + \gamma \mathfrak{W}_{x^\star}\big(h, \mathrm{P}_{T^\perp_{x^\star}} \nabla \widetilde{\Psi}(x^\star)\big) \\ &= \gamma \mathrm{P}_{T_{x^\star}} \nabla^2 \widetilde{R}(x^\star) h + \gamma \mathfrak{W}_{x^\star}\big(h, \mathrm{P}_{T^\perp_{x^\star}} \nabla \widetilde{\Phi}(x^\star)\big) \\ &= \gamma \nabla^2_{\mathcal{M}_{x^\star}} \Phi(x^\star) \mathrm{P}_{T_{x^\star}} h - Hh = Uh, \end{aligned}$$

where $\widetilde{\cdot}$ is the smooth representative of the corresponding function.
Since $-\nabla F(x^\star) \in \mathrm{ri}(\partial R(x^\star))$, we have from [38, Corollary 5.4] that

$$\partial^2 R\big(x^\star | -\nabla F(x^\star)\big) h = \begin{cases} \nabla^2_{\mathcal{M}_{x^\star}} \Psi(x^\star) h + T^\perp_{x^\star}, & h \in T_{x^\star}, \\ \emptyset, & h \notin T_{x^\star}, \end{cases}$$

where $\partial^2 R(x^\star | -\nabla F(x^\star))$ denotes the Mordukhovich generalized Hessian mapping of function $R$ at $(x^\star, -\nabla F(x^\star)) \in \mathrm{gph}\,(\partial R)$ [43]. As $R \in \Gamma_0(\mathbb{R}^n)$, $\partial R$ is a maximal monotone operator, and in view of [50, Theorem 2.1] we have that the mapping $\partial^2 R(x^\star | -\nabla F(x^\star))$ is positive semi-definite, whence we conclude that $\forall h \in T_{x^\star}$,

$$0 \le \gamma\langle \partial^2 R\big(x^\star | -\nabla F(x^\star)\big)h,\, h\rangle = \gamma\langle \nabla^2_{\mathcal{M}_{x^\star}} \Psi(x^\star)h,\, h\rangle = \langle Uh,\, h\rangle.$$

(ii) In this case, $U = \gamma \mathrm{P}_{T_{x^\star}} \nabla^2 \widetilde{R}(x^\star) \mathrm{P}_{T_{x^\star}}$. Let $x_t = x^\star + th$, $t > 0$, for any scalar $t$ and $h \in T_{x^\star}$. Obviously, $x_t \in x^\star + T_{x^\star} = \mathcal{M}_{x^\star}$, and for $t$ sufficiently small, by Fact 3.2, $T_{x_t} = T_{x^\star}$. Thus, $\forall u \in \partial R(x^\star)$ and $\forall v \in \partial R(x_t)$

$$\begin{aligned}
0 \le t^{-2}\langle v - u,\, x_t - x^\star\rangle &= t^{-1}\langle v - u,\, \mathrm{P}_{T_{x^\star}} h\rangle \\
&= t^{-1}\langle \mathrm{P}_{T_{x^\star}}(v - u),\, h\rangle \\
&= t^{-1}\langle \mathrm{P}_{T_{x_t}} v - \mathrm{P}_{T_{x^\star}} u,\, h\rangle \\
\text{(by Fact 3.3)} \quad &= \langle t^{-1}(\nabla_{\mathcal{M}_{x^\star}} R(x_t) - \nabla_{\mathcal{M}_{x^\star}} R(x^\star)),\, h\rangle \\
\text{(by (B.2))} \quad &= \langle t^{-1}\mathrm{P}_{T_{x^\star}}(\nabla \widetilde{R}(x^\star + t\mathrm{P}_{T_{x^\star}} h) - \nabla \widetilde{R}(x^\star)),\, h\rangle.
\end{aligned}$$

Since $\widetilde{R}$ is $C^2$, passing to the limit as $t \to 0$ leads to the desired result. ☐

**Proof of Lemma 4.4.**

(i) (a) is proved using the assumptions and Rademacher theorem. (b) and (c) follow from simple linear algebra arguments.

(ii) From Lemma 4.3, we have $WG = W^{1/2}W^{1/2}GW^{1/2}W^{-1/2}$, meaning that $WG$ is similar to $W^{1/2}GW^{1/2}$. The latter is symmetric and obeys

$$\|W^{1/2}GW^{1/2}\| \le \|W^{1/2}\|\|G\|\|W^{1/2}\| < 1,$$

where we used (i)-(b) to get the last inequality. Thus $W^{1/2}GW^{1/2}$ has real eigenvalues in $]-1, 1[$, and so does $WG$ by similarity. The last statement follows using (i)-(c). ☐

We define the iteration-dependent versions of the matrices in (4.2), *i.e.*

$$\begin{aligned}
&H_k = \gamma_k \mathrm{P}_{T_{x^\star}} \nabla^2 F(x^\star) \mathrm{P}_{T_{x^\star}},\ G_k = \mathrm{Id} - H_k,\ U_k = \gamma_k \nabla^2_{\mathcal{M}_{x^\star}} \Phi(x^\star) \mathrm{P}_{T_{x^\star}} - H_k, \\
&M_{k,1} = \big[(1+b)W(G_k - G), -bW(G_k - G)\big], \\
&M_{k,2} = \big[\big((a_k - b_k) - (a - b)\big)W + (b_k - b)WG_k, -\big((a_k - b_k) - (a - b)\big)W - (b_k - b)WG_k\big].
\end{aligned} \tag{B.5}$$

After the finite identification of $\mathcal{M}_{x^\star}$, we have $x_k \in \mathcal{M}_{x^\star}$ for $x_k$ close enough to $x^\star$. Let $T_{x_k}$ be their corresponding tangent spaces, and define $\tau_k : T_{x^\star} \to T_{x_k}$ the parallel translation along the unique geodesic joining from $x_k$ to $x^\star$.

Before proving Proposition 4.5, we first establish the following intermediate result which provides useful estimates.

**Proposition B.3.** *Under the assumptions of Proposition 4.5, we have*

$$\begin{gathered}
\|y_{a,k} - x^\star\| = O(\|d_k\|),\ \|y_{b,k} - x^\star\| = O(\|d_k\|),\ \|r_{k+1}\| = O(\|d_k\|), \\
(\tau_{k+1}^{-1}\mathrm{P}_{T_{x_{k+1}}} - \mathrm{P}_{T_{x^\star}})\big(\nabla F(y_{b,k}) - \nabla F(x_{k+1})\big) = o(\|d_k\|).
\end{gathered} \tag{B.6}$$

*and*

$$\|W(U_k - U)r_{k+1}\| = o(\|d_k\|),\ \|M_{k,1}d_k\| = o(\|d_k\|)\ \textit{and}\ \|M_{k,2}d_k\| = o(\|d_k\|). \tag{B.7}$$

**Proof.** We have

$$\begin{aligned}\|y_{a,k}-x^\star\| &= \|(1+a_k)r_k - a_k r_{k-1}\| \le (1+a_k)\|r_k\| + a_k\|r_{k-1}\| \\ &\le (1+a_k)(\|r_k\|+\|r_{k-1}\|) \le \sqrt{2}(1+a_k)\|d_k\|,\end{aligned} \tag{B.8}$$

whence we get the first and second estimates. In turn, we obtain

$$\begin{aligned}\|r_{k+1}\| &= \|\mathrm{prox}_{\gamma_k R}\big(y_{a,k}-\gamma_k\nabla F(y_{b,k})\big) - \mathrm{prox}_{\gamma_k R}\big(x^\star-\gamma_k\nabla F(x^\star)\big)\| \\ &\le \|(y_{a,k}-x^\star) - \gamma_k\big(\nabla F(y_{b,k})-\nabla F(x^\star)\big)\| \\ &\le \|(1+a_k)r_k - a_k r_{k-1}\| + \tfrac{\gamma_k}{\beta}\|(1+b_k)r_k - b_k r_{k-1}\| \\ &\le (1+a_k)\|r_k\| + a_k\|r_{k-1}\| + (1+b_k)\tfrac{\gamma_k}{\beta}\|r_k\| + \frac{b_k\gamma_k}{\beta}\|r_{k-1}\| \\ &\le \big((1+a_k)+(1+b_k)\tfrac{\gamma_k}{\beta}\big)(\|r_k\|+\|r_{k-1}\|) \\ &\le \big((1+a_k)+(1+b_k)\tfrac{\gamma_k}{\beta}\big)\sqrt{2}\|d_k\|,\end{aligned} \tag{B.9}$$

where we used non-expansiveness of the proximity operator and assumption (**H.2**). This yields the third estimate. Combining Lemma B.1, assumption (**H.2**), (B.8) and (B.9), we get

$$\begin{aligned}(\tau_{k+1}^{-1}\mathrm{P}_{T_{x_{k+1}}} - \mathrm{P}_{T_{x^\star}})\big(\nabla F(y_{b,k}) - \nabla F(x_{k+1})\big) &= o(\|\nabla F(y_{b,k}) - \nabla F(x_{k+1})\|) \\ &= o(\|y_{b,k}-x^\star\|) + o(\|r_{k+1}\|) = o(\|d_k\|).\end{aligned}$$

Let's now turn to (B.7). Recall the function $\Psi$ defined in the proof of Lemma 4.3(i). First, we have

$$\begin{aligned}\lim_{k\to\infty}\frac{\|W(U_k-U)r_{k+1}\|}{\|r_{k+1}\|} &= \lim_{k\to\infty}\frac{\|W(\gamma_k-\gamma)\nabla^2_{\mathcal{M}_{x^\star}}\Psi(x^\star)\mathrm{P}_{T_{x^\star}}r_{k+1}\|}{\|r_{k+1}\|} \\ &\le \lim_{k\to\infty}|\gamma_k-\gamma|\|W\|\|\nabla^2_{\mathcal{M}_{x^\star}}\Psi(x^\star)\mathrm{P}_{T_{x^\star}}\| = 0,\end{aligned}$$

which entails $\|W(U_k-U)r_{k+1}\| = o(\|r_{k+1}\|) = o(\|d_k\|)$. Again, since $\gamma_k \to \gamma$,

$$\begin{aligned}\lim_{k\to\infty}\frac{\|M_{k,1}d_k\|}{\|d_k\|} &= \lim_{k\to\infty}\frac{\|(1+b)W(G_k-G)r_k - bW(G_k-G)r_{k-1}\|}{\|d_k\|} \\ &\le \lim_{k\to\infty}\frac{(1+b)\|W\|\|G_k-G\|(\|r_k\|+\|r_{k-1}\|)}{\|d_k\|} \\ &\le \lim_{k\to\infty}\frac{(1+b)\|W\||\gamma_k-\gamma|\|\mathrm{P}_{T_{x^\star}}\nabla^2F(x^\star)\mathrm{P}_{T_{x^\star}}\|\sqrt{2}\|d_k\|}{\|d_k\|} \\ &= \lim_{k\to\infty}\sqrt{2}|\gamma_k-\gamma|\big((1+b)\|W\|\|\mathrm{P}_{T_{x^\star}}\nabla^2F(x^\star)\mathrm{P}_{T_{x^\star}}\|\big) = 0,\end{aligned}$$

as $(1+b)\|W\|\|\mathrm{P}_{T_{x^\star}}\nabla^2F(x^\star)\mathrm{P}_{T_{x^\star}}\|$ is obviously bounded (by $2/\beta$). Similarly, for $M_{k,2}$, since $a_k \to a$, $b_k \to$

$b$,

$$
\begin{aligned}
\lim_{k\to\infty}\frac{\|M_{k,2}d_k\|}{\|d_k\|} &= \lim_{k\to\infty}\frac{\|\big(((a_k-b_k)-(a-b))W_k+(b_k-b)W_kG_k\big)(r_k-r_{k-1})\|}{\|d_k\|}\\
&\le \lim_{k\to\infty}\frac{(|a_k-a|+|b_k-b|)\|(W_k+W_kG_k)(r_k-r_{k-1})\|}{\|d_k\|}\\
&\le \lim_{k\to\infty}\frac{(|a_k-a|+|b_k-b|)\|W_k(\mathrm{Id}+G_k)\|\|r_k-r_{k-1}\|}{\|d_k\|}\\
&\le \lim_{k\to\infty}\frac{(|a_k-a|+|b_k-b|)\|W_k(\mathrm{Id}+G_k)\|\sqrt{2}\|d_k\|}{\|d_k\|}\\
&= \lim_{k\to\infty}\sqrt{2}(|a_k-a|+|b_k-b|)\|W_k(\mathrm{Id}+G_k)\| = 0,
\end{aligned}
$$

where $W_k, G_k$ are bounded. □

**Proof of Proposition 4.5.** (1.3) and the first-order optimality condition for problem ($\mathcal{P}_{\mathrm{opt}}$) are respectively equivalent to

$$
\begin{aligned}
y_{a,k}-x_{k+1}-\gamma_k\big(\nabla F(y_{b,k})-\nabla F(x_{k+1})\big) &\in \gamma_k\partial\Phi(x_{k+1})\\
0 &\in \gamma_k\partial\Phi(x^\star).
\end{aligned}
$$

Projecting into $T_{x_{k+1}}$ and $T_{x^\star}$, respectively, and using Fact 3.3, leads to

$$
\begin{aligned}
\gamma_k\tau_{k+1}^{-1}\nabla_{\mathcal{M}_{x^\star}}\Phi(x_{k+1}) &= \tau_{k+1}^{-1}\mathrm{P}_{T_{x_{k+1}}}\big(y_{a,k}-x_{k+1}-\gamma_k\big(\nabla F(y_{b,k})-\nabla F(x_{k+1})\big)\big)\\
\gamma_k\nabla_{\mathcal{M}_{x^\star}}\Phi(x^\star) &= 0.
\end{aligned}
$$

Adding both identities, and subtracting $\tau_{k+1}^{-1}\mathrm{P}_{T_{x_{k+1}}}x^\star$ on both sides, we arrive at

$$
\begin{aligned}
&\tau_{k+1}^{-1}\mathrm{P}_{T_{x_{k+1}}}r_{k+1}+\gamma_k\big(\tau_{k+1}^{-1}\nabla_{\mathcal{M}_{x^\star}}\Phi(x_{k+1})-\nabla_{\mathcal{M}_{x^\star}}\Phi(x^\star)\big)\\
&=\tau_{k+1}^{-1}\mathrm{P}_{T_{x_{k+1}}}(y_{a,k}-x^\star)-\gamma_k\tau_{k+1}^{-1}\mathrm{P}_{T_{x_{k+1}}}\big(\nabla F(y_{b,k})-\nabla F(x_{k+1})\big).
\end{aligned}
\tag{B.10}
$$

By virtue of Lemma B.1, we get

$$
\tau_{k+1}^{-1}\mathrm{P}_{T_{x_{k+1}}}r_{k+1} = \mathrm{P}_{T_{x^\star}}r_{k+1}+(\tau_{k+1}^{-1}\mathrm{P}_{T_{x_{k+1}}}-\mathrm{P}_{T_{x^\star}})r_{k+1} = \mathrm{P}_{T_{x^\star}}r_{k+1}+o(\|r_{k+1}\|).
$$

Using [39, Lemma 5.1], we also have

$$
r_{k+1} = \mathrm{P}_{T_{x^\star}}r_{k+1}+o(\|r_{k+1}\|),
$$

and thus

$$
\tau_{k+1}^{-1}\mathrm{P}_{T_{x_{k+1}}}r_{k+1} = r_{k+1}+o(\|r_{k+1}\|) = r_{k+1}+o(\|d_k\|),
\tag{B.11}
$$

where we also used (B.6). Similarly

$$
\begin{aligned}
\tau_{k+1}^{-1}\mathrm{P}_{T_{x_{k+1}}}(y_{a,k}-x^\star) &= \mathrm{P}_{T_{x^\star}}(y_{a,k}-x^\star)+(\tau_{k+1}^{-1}\mathrm{P}_{T_{x_{k+1}}}-\mathrm{P}_{T_{x^\star}})(y_{a,k}-x^\star)\\
&= \mathrm{P}_{T_{x^\star}}(y_{a,k}-x^\star)+o(\|y_{a,k}-x^\star\|)\\
&= \mathrm{P}_{T_{x^\star}}(y_{a,k}-x^\star)+o(\|d_k\|)\\
&= (1+a_k)\mathrm{P}_{T_{x^\star}}r_k-a_k\mathrm{P}_{T_{x^\star}}r_{k-1}+o(\|d_k\|)\\
&= (1+a_k)r_k-a_kr_{k-1}+o(\|r_k\|)+o(\|r_{k-1}\|)+o(\|d_k\|)\\
&= (y_{a,k}-x^\star)+o(\|d_k\|).
\end{aligned}
\tag{B.12}
$$

Moreover owing to Lemma B.2 and (B.6),

$$\begin{aligned}\tau^{-1}\nabla_{\mathcal{M}_{x^\star}}\Phi(x_{k+1}) - \nabla_{\mathcal{M}_{x^\star}}\Phi(x^\star) &= \nabla^2_{\mathcal{M}_{x^\star}}\Phi(x^\star)\mathrm{P}_{T_{x^\star}} r_{k+1} + o(\|r_{k+1}\|) \\ &= \nabla^2_{\mathcal{M}_{x^\star}}\Phi(x^\star)\mathrm{P}_{T_{x^\star}} r_{k+1} + o(\|d_k\|).\end{aligned} \tag{B.13}$$

Therefore, inserting (B.11), (B.12) and (B.13) into (B.10), we obtain

$$\big(\mathrm{Id} + \gamma_k\nabla^2_{\mathcal{M}_{x^\star}}\Phi(x^\star)\mathrm{P}_{T_{x^\star}}\big)r_{k+1} = (y_{a,k} - x^\star) - \gamma_k\tau^{-1}_{k+1}\mathrm{P}_{T_{x_{k+1}}}\big(\nabla F(y_{b,k}) - \nabla F(x_{k+1})\big) + o(\|d_k\|). \tag{B.14}$$

Owing to (B.6) and local $C^2$-smoothness of $F$, we have

$$\begin{aligned}&\tau^{-1}_{k+1}\mathrm{P}_{T_{x_{k+1}}}\big(\nabla F(y_{b,k}) - \nabla F(x_{k+1})\big)\\ &= \mathrm{P}_{T_{x^\star}}\big(\nabla F(y_{b,k}) - \nabla F(x_{k+1})\big) + o(\|d_k\|)\\ &= \mathrm{P}_{T_{x^\star}}\big(\nabla F(y_{b,k}) - \nabla F(x^\star)\big) - \mathrm{P}_{T_{x^\star}}\big(\nabla F(x_{k+1}) - \nabla F(x^\star)\big) + o(\|d_k\|)\\ &= \mathrm{P}_{T_{x^\star}}\nabla^2 F(x^\star)(y_{b,k} - x^\star) + o(\|y_{b,k} - x^\star\|) - \mathrm{P}_{T_{x^\star}}\nabla^2 F(x^\star)r_{k+1} + o(\|r_{k+1}\|) + o(\|d_k\|)\\ &= \mathrm{P}_{T_{x^\star}}\nabla^2 F(x^\star)\mathrm{P}_{T_{x^\star}}(y_{b,k} - x^\star) - \mathrm{P}_{T_{x^\star}}\nabla^2 F(x^\star)\mathrm{P}_{T_{x^\star}}(x_{k+1} - x^\star) + o(\|d_k\|).\end{aligned} \tag{B.15}$$

Injecting (B.15) in (B.14), we get

$$\begin{aligned}&\big(\mathrm{Id} + \gamma_k\nabla^2_{\mathcal{M}_{x^\star}}\Phi(x^\star)\mathrm{P}_{T_{x^\star}} - \gamma_k\mathrm{P}_{T_{x^\star}}\nabla^2 F(x^\star)\mathrm{P}_{T_{x^\star}}\big)r_{k+1}\\ &= (\mathrm{Id} + U_k)r_{k+1} = (y_{a,k} - x^\star) - H_k(y_{b,k} - x^\star) + o(\|d_k\|),\end{aligned} \tag{B.16}$$

which can be further written as,

$$\begin{aligned}(\mathrm{Id} + U_k)r_{k+1} &= (\mathrm{Id} + U)r_{k+1} + (U_k - U)r_{k+1} = (y_{a,k} - x^\star) - H_k(y_{b,k} - x^\star) + o(\|d_k\|)\\ &= \big((1+a_k)r_k - a_k r_{k-1}\big) - H_k\big((1+b_k)r_k - b_k r_{k-1}\big) + o(\|d_k\|)\\ &= \big((1+a_k)r_k - (1+b_k)H_k r_k\big) - \big(a_k r_{k-1} - b_k H_k r_{k-1}\big) + o(\|d_k\|)\\ &= \big((a_k - b_k)\mathrm{Id} + (1+b_k)G_k\big)r_k - \big((a_k - b_k)\mathrm{Id} + b_k G_k\big)r_{k-1} + o(\|d_k\|)\\ &= \big[(a_k - b_k)\mathrm{Id} + (1+b_k)G_k \quad -\big((a_k - b_k)\mathrm{Id} + b_k G_k\big)\big]\, d_k + o(\|d_k\|).\end{aligned}$$

Inverting $\mathrm{Id} + U$ (which is possible thanks to Lemma 4.3), we obtain

$$r_{k+1} + W(U_k - U)r_{k+1} = \big[(a_k - b_k)W + (1+b_k)WG_k \quad -(a_k - b_k)W - b_k WG_k\big]\, d_k + o(\|d_k\|).$$

Using the estimates (B.7), we get

$$\begin{aligned}d_{k+1} &= \begin{bmatrix}(a_k - b_k)W + (1+b_k)WG_k & -(a_k - b_k)W - b_k WG_k\\ \mathrm{Id} & 0\end{bmatrix} d_k + o(\|d_k\|)\\ &= \left(M + \begin{bmatrix}M_{k,1}\\ 0\end{bmatrix} + \begin{bmatrix}M_{k,2}\\ 0\end{bmatrix}\right) d_k + o(\|d_k\|) = Md_k + o(\|d_k\|).\end{aligned}$$

□

**Proof of Proposition 4.7.**

(i) We have

$$\begin{aligned}M\begin{pmatrix}r_1\\ r_2\end{pmatrix} &= \begin{bmatrix}(a-b)\mathrm{Id} + (1+b)G, & -(a-b)\mathrm{Id} - bG\\ \mathrm{Id}, & 0\end{bmatrix}\begin{pmatrix}r_1\\ r_2\end{pmatrix}\\ &= \begin{pmatrix}(a-b)r_1 + (1+b)Gr_1 - (a-b)r_2 - bGr_2\\ r_1\end{pmatrix} = \sigma\begin{pmatrix}r_1\\ r_2\end{pmatrix},\end{aligned}$$

and thus $r_1 = \sigma r_2$. Inserting this in the first identity, we obtain

$$\sigma^2 r_2 = (a-b)\sigma r_2 + (1+b)\sigma G r_2 - (a-b) r_2 - bGr_2$$
$$\iff Gr_2 = \frac{(a-b)(1-\sigma)+\sigma^2}{(1+b)\sigma - b} r_2 = \eta r_2 \implies 0 = \sigma^2 - \big((a-b)+(1+b)\eta\big)\sigma + (a-b) + b\eta.$$

(ii) For this quadratic equation of $\sigma$, the two roots are

$$\sigma_1 = \frac{\big((a-b)+(1+b)\eta\big)+\sqrt{\Delta_\sigma}}{2}, \quad \sigma_2 = \frac{\big((a-b)+(1+b)\eta\big)-\sqrt{\Delta_\sigma}}{2}. \tag{B.17}$$

where $\Delta_\sigma$ is the discriminant

$$\Delta_\sigma = \big((a-b)+(1+b)\eta\big)^2 - 4\big((a-b)+b\eta\big),$$

which is a quadratic function of 3 variables. Consider the following 3 linear functions of $a$

$$\begin{aligned} a_1 &= (1-\eta)b - \eta, \\ a_2 &= (1-\eta)b + (1-\sqrt{1-\eta})^2 \begin{cases} \Delta_\sigma \le 0 : a_2 \le a \le 1 \le (1-\eta)b + (1+\sqrt{1-\eta})^2, \\ \Delta_\sigma \ge 0 : a \le a_2, \end{cases} \\ a_3 &= (1-\eta)b - \frac{1+\eta}{2}. \end{aligned} \tag{B.18}$$

Recall from Lemma 4.4(i) that $\eta \in ]-1,1[$. Thus, $a_1 \ge a_2$ when $\eta \in ]-1,0]$, $a_1 \le a_2$ when $\eta \in [0,1[$, and $a_3$ is smaller than both $a_1, a_2$ independently of $\eta$. We now discuss each case.

**Case** $\eta \in ]-1,0]$**:** We have $a_1 \ge a_2$,

- **Subcase** $a \in [a_2, 1[$**:** $\sigma_{1,2}$ are complex, hence

$$|\sigma|^2 = \frac{\big((a-b)+(1+b)\eta\big)^2 - \Big(\big((a-b)+(1+b)\eta\big)^2 - 4\big((a-b)+b\eta\big)\Big)}{4} = a - b + b\eta. \tag{B.19}$$

  Since $a_2 \le 1 \iff b \le \frac{1-(1-\sqrt{1-\eta})^2}{1-\eta}$, then we have $(1-\sqrt{1-\eta})^2 \le |\sigma|^2 \le 1 + (\eta-1)b < 1$.
- **Subcase** $a \in [0, a_2]$**:** $\Delta_\sigma \ge 0$ and $\sigma_2$ has the bigger absolute value, then

$$\begin{aligned} |\sigma_2| < 1 &\iff -\big((a-b)+(1+b)\eta\big) + \sqrt{\Delta_\sigma} < 2 \\ &\iff \Delta_\sigma < 4 + 4\big((a-b)+(1+b)\eta\big) + \big((a-b)+(1+b)\eta\big)^2 \\ &\iff \frac{2(b-a)-1}{1+2b} < \eta, \end{aligned} \tag{B.20}$$

  which means that $|\sigma_2| \le 1$ for $a \in [a_3, a_2]$, and $|\sigma_2| \ge 1$ for $a \in [0, a_3]$. Moreover, $a_3 \le 0$ for $b \in [0, \frac{1+\eta}{2(1-\eta)}]$, meaning that if $\eta \ge \frac{1}{3}$, then $|\sigma_2| \le 1$ for $a \in [0, a_2]$.

**Case** $\eta \in [0,1[$**:** First we have $a_2 \ge a_1$, and moreover

$$a_1 = 0 \iff b = \frac{\eta}{1-\eta} \begin{cases} \le 1 : \eta \in [0, 0.5], \\ \ge 1 : \eta \in [0.5, 1[. \end{cases}$$

Obviously, we have $|\sigma| \le 1$ holds for any $a \in [0, a_2]$ as long as $\eta \in [0.5, 1]$, though this situation is useless as $b \in [0,1]$. In the subcases hereafter, we only consider $\eta \in [0, 0.5]$.

- **Subcase** $a \in [a_2, 1[$**:** same result as (B.19).
- **Subcase** $a \in [a_1, a_2]$**:** $\sigma_1 \geq |\sigma_2|$, hence

$$\begin{aligned}\sigma_1 < 1 &\iff \big((a-b)+(1+b)\eta\big)+\sqrt{\Delta_\sigma} < 2 \\ &\iff \Delta_\sigma < 4 - 4\big((a-b)+(1+b)\eta\big)+\big((a-b)+(1+b)\eta\big)^2 \\ &\iff 0 < 4(1-\eta).\end{aligned} \tag{B.21}$$

- **Subcase** $a \in [0, a_1]$**:** We have $|\sigma_2| \geq |\sigma_1|$, hence (B.20) applies and the result follows.

Summarizing this discussion yields the claimed result. □

**Proof of Theorem 4.15.** Since $R$ is locally polyhedral, then $\nabla_{\mathcal{M}_{x^\star}} \Phi(x_k)$ is locally constant along $\mathcal{M}_{x^\star} = x^\star + T_{x^\star}$ around $x^\star$. Thus, embarking from (B.16) in the proof of Proposition 4.5, for $k$ large enough, we get

$$x_{k+1} - x^\star = (y_{a,k} - x^\star) - E_k(y_{b,k} - x^\star),$$

where we used the mean-value theorem with $E_k = \gamma_k \int_0^1 \nabla^2 F(x^\star + t(y_{b,k} - x^\star))\mathrm{d}t \succeq 0$. Using that $E_k$ is symmetric and $\mathrm{Im}(E_k)^\perp = V$, we have

$$\mathrm{P}_V(x_{k+1} - x^\star) = \mathrm{P}_V(y_{a,k} - x^\star) = (1+a_k)\mathrm{P}_V(x_k - x^\star) - a_k \mathrm{P}_V(x_{k-1} - x^\star).$$

If $a_k = 0$, then $\mathrm{P}_V(x_{k+1} - x^\star) = \mathrm{P}_V(x_k - x^\star)$. Thus, in the rest, without loss of generality, we assume that $a_k > 0$ for $k$ large enough. The above iteration leads to

$$\begin{pmatrix}\mathrm{P}_V(x_{k+1} - x^\star) \\ \mathrm{P}_V(x_k - x^\star)\end{pmatrix} = \begin{bmatrix}(1+a_k)\mathrm{Id} & -a_k\mathrm{Id} \\ \mathrm{Id} & 0\end{bmatrix}\begin{pmatrix}\mathrm{P}_V(x_k - x^\star) \\ \mathrm{P}_V(x_{k-1} - x^\star)\end{pmatrix}.$$

It is straightforward to check that $N_k \stackrel{\text{def}}{=} \begin{bmatrix}(1+a_k)\mathrm{Id} & -a_k\mathrm{Id} \\ \mathrm{Id} & 0_n\end{bmatrix}$ is invertible and admits two eigenvalues $a_k > 0$ and $1$ respectively. Iterating the above argument, and owing to the fact that $x_k, y_{a,k}, y_{b,k} \to x^\star$, we get

$$\begin{pmatrix}0 \\ 0\end{pmatrix} = \left(\prod_{j=k}^{\infty} N_j\right)\begin{pmatrix}\mathrm{P}_V(x_k - x^\star) \\ \mathrm{P}_V(x_{k-1} - x^\star)\end{pmatrix},$$

and $\prod_{j=k}^{\infty} N_j$ is invertible. Therefore, we obtain that $x_k - x^\star \in V^\perp$, and in turn, $y_{a,k} - x^\star \in V^\perp$ and $y_{b,k} - x^\star \in V^\perp$, for all large enough $k$. Observe that $V^\perp \subset T_{x^\star}$, it then follows that

$$x_{k+1} - x^\star = y_{a,k} - x^\star - \mathrm{P}_{V^\perp} E_k \mathrm{P}_{V^\perp}(y_{b,k} - x^\star).$$

By definition, $\mathrm{P}_{V^\perp} E_k \mathrm{P}_{V^\perp}$ is symmetric positive definite. Thus, replacing $H_k$ by $\mathrm{P}_{V^\perp} E_k \mathrm{P}_{V^\perp}$, $G$ and $M$ accordingly, in Lemma 4.4 and Corollary 4.9, and applying Theorem 4.11 leads to the result. □

# References

[1] P-A. Absil, R. Mahony, and R. Sepulchre. *Optimization algorithms on matrix manifolds*. Princeton University Press, 2009.

[2] P-A. Absil, R. Mahony, and J. Trumpf. An extrinsic look at the Riemannian Hessian. In *Geometric Science of Information*, pages 361–368. Springer, 2013.

[3] A. Agarwal, S. Negahban, and M. J. Wainwright. Fast global convergence of gradient methods for high-dimensional statistical recovery. *The Annals of Statistics*, 40(5):2452–2482, 2012.

[4] F. Alvarez. On the minimizing property of a second order dissipative system in Hilbert spaces. *SIAM Journal on Control and Optimization*, 38(4):1102–1119, 2000.

[5] F. Alvarez and H. Attouch. An inertial proximal method for maximal monotone operators via discretization of a nonlinear oscillator with damping. *Set-Valued Analysis*, 9(1-2):3–11, 2001.

[6] H. Attouch, Z. Chbani, J. Peypouquet, and P. Redont. Fast convergence of inertial dynamics and algorithms with asymptotic vanishing damping. Technical Report Optimization online 5179, 2015.

[7] H. Attouch and J. Peypouquet. The rate of convergence of Nesterov's accelerated Forward–Backward method is actually $o(k^{-2})$. Technical Report arXiv:1510.08740, 2015.

[8] H. Attouch, J. Peypouquet, and P. Redont. A dynamical approach to an inertial Forward–Backward algorithm for convex minimization. *SIAM J. Optim.*, 24(1):232–256, 2014.

[9] H. Attouch, J. Peypouquet, and P. Redont. On the fast convergence of an inertial gradient-like dynamics with vanishing viscosity. Technical Report arXiv:1507.04782, 2015.

[10] J.-F. Aujol and C. Dossal. Stability of over-relaxations for the Forward–Backward algorithm, application to fista. *SIAM Journal on Optimization*, 25(4):2408–2433, 2015.

[11] J. B. Baillon and G. Haddad. Quelques propriétés des opérateurs angle-bornés etn-cycliquement monotones. *Israel Journal of Mathematics*, 26(2):137–150, 1977.

[12] H. Bauschke and P. L. Combettes. *Convex Analysis and Monotone Operator Theory in Hilbert Spaces*. Springer, 2011.

[13] A. Beck and M. Teboulle. Fast gradient-based algorithms for constrained total variation image denoising and deblurring problems. *Image Processing, IEEE Transactions on*, 18(11):2419–2434, 2009.

[14] A. Beck and M. Teboulle. A fast iterative shrinkage-thresholding algorithm for linear inverse problems. *SIAM Journal on Imaging Sciences*, 2(1):183–202, 2009.

[15] J. Bolte, A. Daniilidis, and A. Lewis. The Łojasiewicz inequality for nonsmooth subanalytic functions with applications to subgradient dynamical systems. *SIAM J. Optim*, 17(4):1205–1223, 2006.

[16] K. Bredies and D. A. Lorenz. Linear convergence of iterative soft-thresholding. *Journal of Fourier Analysis and Applications*, 14(5-6):813–837, 2008.

[17] Regina S. Burachik and Alfredo N. Iusem. *Set-valued Mappings and Enlargements of Monotone Operators*. Optimization and Its Applications. Springer, 2008.

[18] E. J. Candès and B. Recht. Simple bounds for recovering low-complexity models. *Mathematical Programming*, 141(1-2):577–589, 2013.

[19] A. Chambolle and C. Dossal. On the convergence of the iterates of the "fast iterative shrinkage/thresholding algorithm". *Journal of Optimization Theory and Applications*, 166(3):968–982, 2015.

[20] I. Chavel. *Riemannian geometry: a modern introduction*, volume 98. Cambridge University Press, 2006.

[21] P. L. Combettes. Quasi-Fejérian analysis of some optimization algorithms. *Studies in Computational Mathematics*, 8:115–152, 2001.

[22] P. L. Combettes and B. C. Vũ. Variable metric Forward–Backward splitting with applications to monotone inclusions in duality. *Optimization*, 63(9):1289–1318, 2014.

[23] A. Daniilidis, D. Drusvyatskiy, and A. S. Lewis. Orthogonal invariance and identifiability. *SIAM Journal on Matrix Analysis and Applications*, 35(2):580–598, 2014.

[24] A. Daniilidis, W. Hare, and J. Malick. Geometrical interpretation of the predictor-corrector type algorithms in structured optimization problems. *Optimization: A Journal of Mathematical Programming & Operations Research*, 55(5-6):482–503, 2009.

[25] V. Duval and G. Peyré. Sparse spikes deconvolution on thin grids. *arXiv preprint arXiv:1503.08577*, 2015.

[26] T. Goldstein, B. O'Donoghue, S. Setzer, and R. Baraniuk. Fast alternating direction optimization methods. *SIAM Journal on Imaging Sciences*, 7(3):1588–1623, 2014.

[27] M. Gu, L.-H. Lim, and C. J. Wu. ParNes: a rapidly convergent algorithm for accurate recovery of sparse and approximately sparse signals. *Numerical Algorithms*, 64(2):321–347, 2012.

[28] E. Hale, W. Yin, and Y. Zhang. Fixed-point continuation for $\ell_1$-minimization: methodology and convergence. *SIAM Journal on Optimization*, 19(3):1107–1130, 2008.

[29] W. L. Hare. Identifying active manifolds in regularization problems. In H. H. Bauschke, R. S., Burachik, P. L. Combettes, V. Elser, D. R. Luke, and H. Wolkowicz, editors, *Fixed-Point Algorithms for Inverse Problems in Science and Engineering*, volume 49 of *Springer Optimization and Its Applications*, chapter 13. Springer, 2011.

[30] W. L. Hare and A. S. Lewis. Identifying active constraints via partial smoothness and prox-regularity. *Journal of Convex Analysis*, 11(2):251–266, 2004.

[31] W. L. Hare and A. S. Lewis. Identifying active manifolds. *Algorithmic Operations Research*, 2(2):75–82, 2007.

[32] J.-B. Hiriart-Urruty and C. Lemaréchal. *Convex Analysis And Minimization Algorithms*, volume I and II. Springer, 2001.

[33] K. Hou, Z. Zhou, A. M.-C. So, and Z.Q. Luo. On the linear convergence of the proximal gradient method for trace norm regularization. In *Advances in Neural Information Processing Systems*, pages 710–718, 2013.

[34] P. R. Johnstone and P. Moulin. A Lyapunov analysis of FISTA with local linear convergence for sparse optimization. *arXiv preprint arXiv:1502.02281*, 2015.

[35] J. M. Lee. *Smooth manifolds*. Springer, 2003.

[36] C. Lemaréchal, F. Oustry, and C. Sagastizábal. The U-Lagrangian of a convex function. *Trans. Amer. Math. Soc.*, 352(2):711–729, 2000.

[37] A. S. Lewis. Active sets, nonsmoothness, and sensitivity. *SIAM Journal on Optimization*, 13(3):702–725, 2003.

[38] A. S. Lewis and S. Zhang. Partial smoothness, tilt stability, and generalized Hessians. *SIAM Journal on Optimization*, 23(1):74–94, 2013.

[39] J. Liang, J. Fadili, and G. Peyré. Local linear convergence of Forward–Backward under partial smoothness. In *Advances in Neural Information Processing Systems*, pages 1970–1978, 2014.

[40] P. L. Lions and B. Mercier. Splitting algorithms for the sum of two nonlinear operators. *SIAM Journal on Numerical Analysis*, 16(6):964–979, 1979.

[41] D. A. Lorenz and T. Pock. An accelerated Forward–Backward algorithm for monotone inclusions. *arXiv preprint arXiv:1403.3522*, 2014.

[42] S. A. Miller and J. Malick. Newton methods for nonsmooth convex minimization: connections among-Lagrangian, Riemannian Newton and SQP methods. *Mathematical programming*, 104(2-3):609–633, 2005.

[43] B.S. Mordukhovich. Sensitivity analysis in nonsmooth optimization. *Theoretical Aspects of Industrial Design (D. A. Field and V. Komkov, eds.), SIAM Volumes in Applied Mathematics*, 58:32–46, 1992.

[44] A. Moudafi and M. Oliny. Convergence of a splitting inertial proximal method for monotone operators. *Journal of Computational and Applied Mathematics*, 155(2):447–454, 2003.

[45] Y. Nesterov. A method for solving the convex programming problem with convergence rate $O(1/k^2)$. *Dokl. Akad. Nauk SSSR*, 269(3):543–547, 1983.

[46] Y. Nesterov. *Introductory lectures on convex optimization: A basic course*, volume 87. Springer, 2004.

[47] Y. Nesterov. Gradient methods for minimizing composite objective function. 2007.

[48] B. O'Donoghue and E. Candes. Adaptive restart for accelerated gradient schemes. *Foundations of computational mathematics*, 15(3):715–732, 2015.

[49] Z. Opial. Weak convergence of the sequence of successive approximations for nonexpansive mappings. *Bulletin of the American Mathematical Society*, 73(4):591–597, 1967.

[50] R. A. Poliquin and R. T. Rockafellar. Tilt stability of a local minimum. *SIAM Journal on Optimization*, 8(2):287–299, 1998.

[51] B. T. Polyack. Some methods of speeding up the convergence of iterative methods. *Zh. Vychisl. Mat. Mat. Fiz.*, 4:1–17, 1964.

[52] B. T. Polyak. *Introduction to optimization*. Optimization Software, 1987.

[53] R. T. Rockafellar. Monotone operators and the proximal point algorithm. *SIAM Journal on Control and Optimization*, 14(5):877–898, 1976.

[54] R. T. Rockafellar and R. Wets. *Variational analysis*, volume 317. Springer Verlag, 1998.

[55] L. I. Rudin, S. Osher, and E. Fatemi. Nonlinear total variation based noise removal algorithms. *Physica D: Nonlinear Phenomena*, 60(1):259–268, 1992.

[56] S. T. Smith. Optimization techniques on Riemannian manifolds. *Fields institute communications*, 3(3):113–135, 1994.

[57] B. F. Svaiter and R. S. Burachik. $\varepsilon$-enlargements of maximal monotone operators in banach spaces. *Set-Valued Anal.*, 7:117–132, 1999.

[58] S. Tao, D. Boley, and S. Zhang. Local linear convergence of ISTA and FISTA on the LASSO problem. *arXiv preprint arXiv:1501.02888*, 2015.

[59] P. Tseng and S. Yun. A coordinate gradient descent method for nonsmooth separable minimization. *Math. Prog. (Ser. B)*, 117, 2009.

[60] S. Vaiter, C. Deledalle, J. M. Fadili, G. Peyré, and C. Dossal. The degrees of freedom of partly smooth regularizers. *Annals of the Institute of Mathematical Statistics*, 2015. to appear.

[61] S. Vaiter, M. Golbabaee, J. Fadili, and G. Peyré. Model selection with low complexity priors. *Information and Inference*, page iav005, 2015.

[62] S. Vaiter, G. Peyré, and J. Fadili. Model consistency of partly smooth regularizers. *arXiv preprint arXiv:1405.1004*, 2014.

[63] S. J. Wright. Identifiable surfaces in constrained optimization. *SIAM Journal on Control and Optimization*, 31(4):1063–1079, 1993.